\documentclass[preprint,11pt]{elsarticle}
\usepackage[toc,page]{appendix}

\usepackage{amsmath,amsthm,amsfonts,bm,hyperref}

\usepackage{pgfplots}
\usetikzlibrary{patterns}
\pgfplotsset{compat=newest}

\usepackage{subcaption}
\usepackage{float,indentfirst,algorithm}

\usepackage{booktabs}
\usepackage[margin=1in]{geometry}

\theoremstyle{plain}

\newtheorem{remark}{Remark}

\usepackage{listings} 
\usepackage[capitalise]{cleveref}
\usepackage{multirow}

\definecolor{BurntOrange}{rgb}{0.8, 0.33, 0.0}
\usepackage{xcolor}

\usepackage{url}

\newcommand{\Mcal}{\ensuremath{\mathcal{M}}}
\newcommand{\Fcal}{\ensuremath{\mathcal{F}}}
\newcommand{\Pcal}{\ensuremath{\mathcal{P}}}
\newcommand{\Tcal}{\ensuremath{\mathcal{T}}}
\newcommand{\Ccal}{\ensuremath{\mathcal{C}}}

\newcommand{\bt}[0]{\bm{\theta}}
\newcommand{\bx}{\mathbf{x}}
\newcommand{\bu}{\mathbf{u}}
\newcommand{\bv}[0]{\mathbf{v}}
\newcommand{\bw}{\mathbf{w}}
\newcommand{\blam}[0]{\bm{\lambda}}
\newcommand{\bSig}[0]{\bm{\Sigma}}
\newcommand{\bP}[0]{\mathbf{P}}
\newcommand{\bF}[0]{\mathbf{F}}
\newcommand{\bK}[0]{\mathbf{K}}
\DeclareMathOperator{\E}{\mathbb{E}}
\newcommand{\RR}[0]{\mathbb{R}}
\DeclareMathOperator{\R}{\mathbb{R}}
\newcommand{\bW}{\mathbf{W}}
\usepackage{tcolorbox}
\usepackage[authormarkup=none,final]{changes}
\definechangesauthor[name=Reviewer 1, color=blue]{r1}
\definechangesauthor[name=Reviewer 2, color=orange]{r2}

\begin{document}

\begin{abstract}
    We present a new approach for predictive modeling and its uncertainty quantification for mechanical systems, where coarse-grained models such as constitutive relations are derived directly from observation data. We explore the use of a neural network to represent the unknown constitutive relations, compare the neural networks with piecewise linear functions, radial basis functions, and radial basis function networks, and show that the neural network outperforms the others in certain cases. We analyze the approximation error of the neural networks using a scaling argument. The training and predicting processes in our framework combine the finite element method, automatic differentiation, and neural networks (or other function approximators).  Our framework also allows uncertainty quantification in the form of confidence intervals. Numerical examples on a multiscale fiber-reinforced plate problem and a nonlinear rubbery membrane problem from solid mechanics demonstrate the effectiveness of our framework.
\end{abstract}


\begin{frontmatter}
\title{Learning Constitutive Relations from Indirect Observations Using Deep Neural Networks
    }
\author[rvt1]{Daniel~Z.~Huang\corref{cor1}}
\ead{zhengyuh@stanford.edu}

\author[rvt1]{Kailai~Xu\corref{cor1}}
\ead{kailaix@stanford.edu}

\author[rvt1,rvt2,rvt3]{Charbel~Farhat}
\ead{cfarhat@stanford.edu}

\author[rvt1,rvt2]{Eric~Darve}
\ead{darve@stanford.edu}

\address[rvt1]{Institute for Computational and Mathematical Engineering,
               Stanford University, Stanford, CA, 94305}
\address[rvt2]{Mechanical Engineering, Stanford University, Stanford, CA, 94305}
\address[rvt3]{Aeronautics and Astronautics, Stanford University, Stanford, CA, 94305}
\cortext[cor1]{Both authors contributed equally to this work, and are listed in alphabetical order.}

\end{frontmatter}

\section{Introduction}

Many practical problems arising from various engineering and scientific applications are heterogeneous and multi-scale in nature. The simulations of such problems based on first principles remain prohibitively expensive. Coarse-grained models are often applied to approximate the effect of the microscopic interactions, to simplify and accelerate these simulations. For example, in solid mechanics, the constitutive relations might be derived from interaction forces between atoms. The constitutive relations can also be modeled empirically based on theoretical knowledge with ideal assumptions and calibrated on limited tensile test data in coarse scales. These modeling efforts can lead to affordable simulations of large scale engineering and scientific applications.

There are two kinds of coarse-grained models in general: (1) \textit{purely phenomenological models}, which directly relate several different empirical observations of phenomena to each other; (2)\textit{ multi-scale models}, which are derived from finer scales or even from the molecular-based theories. The present work belongs to the former one; traditionally it leads to certain model forms with a few parameters (e.g., Young's modulus) to be estimated with inverse analysis~\cite{trujillo1997practical}. We focus on calibrating the coarse-grained model, specifically the constitutive relations in solid mechanics. There are in general two ways to calibrate the constitutive relations:

\begin{enumerate}
    \item Methods relying on \textit{direct} data, such as strain-stress pairs, strain-strain energy pairs, or strain-stress increments pairs. 
    
    For example, neural networks have been used to model the constitutive relations in a variety of materials, including concrete~\cite{ghaboussi1991knowledge}, sands~\cite{ellis1995stress}, hyper-elastic materials~\cite{shen2005finite}, nonlinear elastic composites~\cite{le2015computational}, crystal elastic materials~\cite{ling2016machine}, viscoelastic materials~\cite{furukawa1998implicit} and even multi-scale porous materials~\cite{wang2018multiscale}.
    
    These data points consist of experimental measurements and numerical simulation results, which are generated from sub-scale simulations, like representative volume element (RVE) simulations~\cite{hashin1983analysis, sun1996prediction, feyel2000fe2, kanit2003determination, yuan2008toward}, or post-processed from direct numerical simulations which need to resolve all scales of the problem. 
    
    However, the comprehensive strain-stress relation measurement relying on simple mechanical tests, such as tensile or bending tests, is challenging, especially for anisotropic materials. For the RVE approach, which can efficiently generate direct data to train neural networks~\cite{le2015computational, bessa2017framework, ling2016machine, wang2018multiscale}, the determination of the RVE size~\cite{kanit2003determination} and finer-scale material properties may be difficult. Direct numerical simulations generate high-fidelity training data sets, but for most practical problems, their computational costs are still unaffordable.

    \item Methods relying on \textit{indirect} data, such as deformations of structures under different load conditions.
    
    Deformations are measured by techniques such as digital image correlation or grid method~\cite{surrel1994moire}.
    These techniques can record complete heterogeneous fields, which are rich in the constitutive relations. However these data are {\it indirect}, i.e., there is generally no closed-form solution allowing a direct link between measurements and the stress or the underlying constitutive relations. The virtual fields method~\cite{grediac2006virtual, avril2008overview, geymonat2002identification, feng2001genetic} has been designed to apply the finite element method (FEM) to bridge the full-field data with parametric constitutive relations. An inverse analysis is used to identify these constitutive parameters. \deleted[id=r2]{Tartakovsky et al.~\cite{tartakovsky2018learning} applied deep neural networks to directly inform the unknown constitutive relationship in the non-linear diffusion equations from the full-field data.} 
\end{enumerate}

We propose an approach that combines traditional numerical discretization schemes and data-driven functional approximation for predictive modeling relying on \textit{indirect} data. Our goal is to build a coarse-grained constitutive relation model with the following workflow. We conduct various experiments: given a mechanical system, we apply different boundary conditions such as external forces and observe the resulting deformation field. We wish to learn a constitutive relation model (possibly a coarse-grained or low-fidelity model) that can reproduce the observed deformations. Generally speaking, a constitutive relation is a function $\Mcal(u, \bx)$ where $\bx$ is a location and $u(\bx)$ is a function that represents the deformation of the object.\footnote{More general relations include $\Mcal(u, \dot{u}, \bx)$, etc, but they will be the object of future work.} 
Given some boundary conditions, the deformation $u$, and the function $\Mcal(\bullet, \bullet)$, we can compute the internal forces from mechanics. 

Our training set to learn the constitutive model $\Mcal$ consists of a series of ``experiments'' (we use quotes since in practice these experiments may correspond to numerical calculations using, for example, a fine-scale or even atomistic model) with different choices of boundary conditions. For each ``experiment,'' we record $u(\bx)$ at discrete locations. This problem is therefore not a regression problem where we must fit some given data $(u_i,\bx_i) \to \Mcal_i$. Instead, we consider that we have a functional form $\Mcal_{\bt}(u)$, which is homogenized spatially, or $\Mcal_{\bt}(u, \bx)$, when the underlying model varies significantly in $\bx$. 
\added[id=r1]{Here, $\bt$ denotes the parameters in the data-driven constitutive model.}
However the latter one can be computational expensive since we need to construct a constitutive relation for each $\bx$; consequently, this case is only used for analysis in the present work. 

Using \added[id=r2]{automatic differentiation} and gradient descent methods, the parameters $\bt$ will be optimized such that the internal force $\bP$---derived from  $\Mcal_{\bt}$ and our observations $\bu_i$---matches the imposed external forces $\bF_i$.
\begin{equation*}
    \min_{\bt} \sum_{i}\left(\bP(\bu_i, \Mcal_{\bt}(\bu_i)) - \bF_i \right)^2
\end{equation*}
\added[id=r2]{ This supports a broad range of traditional model calibration problems such as parameter estimation. Moreover, this approach is capable of learning unknown model forms.}


A critical element is to determine the functional form best suited to this task. Many functional forms rely on a partitioning of the space $(u,\bx)$ into cells $\Omega_i$ (for example in simplices, parallelepipeds, or simple geometrical shapes) and using low-order polynomials, or a local Fourier basis inside each $\Omega_i$~\cite{boyd2001chebyshev,canuto2006spectral,canuto2012spectral}. However, if we examine our approach, we realize that these techniques are ill-suited. 

Indeed, consider a particular experiment associated with given boundary conditions. This will lead to some field $u$, from which strain values can be computed. However, these strain values are not ``uniformly'' distributed in the domain. They typically will lie on some low-dimensional manifold (see \Cref{fig:membrane_strain_paris} for example). As a result, the collection of all strain values observed throughout all experiments has a very irregular and ``anisotropic'' distribution (see \Cref{fig:membrane_strain_paris} to see a specific example of what we mean). Building an appropriate $\Omega_i$ is therefore challenging and error-prone. Choosing large $\Omega_i$ (large diameter) leads to poor reconstruction while small $\Omega_i$ may lead to instabilities if no or few sample (``training'') points fall inside the cell. If one uses a Delaunay-type triangulation~\cite{field1988laplacian,lee1980two}, the elements will be extremely distorted leading to ill-conditioned numerical computations. 

The order of the basis (e.g., order of the polynomials) needs to be chosen carefully. A low-order basis leads to large errors while a high-order basis leads to an unstable interpolation procedure.

Although such approaches may provide an accurate reconstruction of $\Mcal$ near observed points, their accuracy typically deteriorates as we move away. For example, Chebyshev polynomials~\cite{berrut2004barycentric,mason2002chebyshev} are known to diverge rapidly outside the $[-1,1]$ interval. In addition, such approaches do not extend well to high-dimensional input data because the basis construction typically relies on a tensor product construction which leads to an exponential number of basis functions in the dimension of the space (that is the number of the basis functions scales like $O(p^d)$ where $p$ is the order and $d$ the dimension).

A more natural choice for such problems is to use radial basis functions~\cite{de2007mesh,park1991universal,rippa1999algorithm,schaback1995error}, that is an approximation of the type:
\[ f(x) \approx \sum_i \alpha_i \; g_\sigma(\|\bx - \bx_i\|) \]
where $g_\sigma$ is, for example, a Gaussian function, an exponential, or a multiquadrics~\cite{mouat2002rbf,genton2015cross}, and $x_i$ are centers used in the approximation; $\sigma$ is a scale parameter (for example the standard deviation of the Gaussian or decay rate of the exponential). Such approaches work quite well even in high-dimension. However, they have many drawbacks. The main one is probably that computing the coefficients $\alpha_i$ requires computing with the matrix $a_{ij} \overset{\text{def}}{=} g_\sigma(\|\bx_i - \bx_j\|)$ which is known to become ill-conditioned as the centers $\bx_i$ get close to each other. As a result, even a small perturbation or error in the input data will lead to large changes in the model coefficients, which is undesirable. Methods like Kriging or Bayesian approaches~\cite{rasmussen2003gaussian,stein2012interpolation,van2004kriging,williams1996gaussian,bernardo1998regression,gibbs1998bayesian,williams1998bayesian} require an {\it a priori} statistical distribution which may or may not apply to the problem at hand. \replaced[id=r1]{In these methods, accounting for observation uncertainties leads to the addition of a positive term (nugget) along the diagonal of the covariance matrix which makes the linear system better conditioned.}{The prior information typically leads to better conditioned linear systems that are easier to solve.} We note however that in Gaussian Process Regression the matrix used to calculate the model reverts to $g_\sigma(\|\bx_i - \bx_j\|)$ in the absence of noise in the observation and, as the density of $\bx_i$ increases, the ill-conditioning appears again.

In this context, deep neural networks~(DNNs) offer many advantages. \deleted[id=r2]{They possess ``universal approximation'' properties}~\cite{lorentz1996constructive, mhaskar1996neural, xu2018calibrating, 2019arXiv190107758X}. They can be trained using nonuniform point cloud data. Using appropriate regularization, they have good ``generalization'' properties, i.e., they remain accurate even away from training points. For example, as the regularization penalization factor increases (using $\ell_2$ or $\ell_1$ regularization), the regression function from a neural network becomes more ``linear'' and flat away from training points. Neural networks are known to work well even for complex, highly inhomogeneous or anisotropic distribution of training points (that is dense along certain directions and sparse along others)~\cite{lorentz1996constructive, mhaskar1996neural, xu2018calibrating, 2019arXiv190107758X}. For example, we will show in our benchmarks that DNNs outperform piecewise linear functions~(\cref{sect:cplrs}), radial basis functions~(\cref{sect:crbf}), and radial basis function networks~(\cref{sect:rbfn}).  \added[id=r2]{It is worth mentioning that Tartakovsky et al.~\cite{tartakovsky2018learning} applied  physics-informed neural networks (PINN) \cite{raissi2019physics} to learn the unknown constitutive relations in the non-linear diffusion equations from the full-field data. However, in their approach both the solution (displacements) and the constitutive relation are approximated by neural networks, while our approach only approximates the constitutive relation with neural networks and then couples it with finite element solvers for the mechanics equations.}

The applicability and accuracy of the learning procedure are analyzed based on a model problem. For problems with a smooth underlying constitutive relation, the learning process delivers an approximate constitutive relation with an error bound that depends on both the optimization error tolerance and the numerical partial differential equation (PDE) discretization error. Moreover, the uncertainties in the constitutive relation are also learned by our algorithm during the training process. Through sensitivity analysis, the uncertainty can be used to provide error bounds and intervals of confidence for the prediction.
\deleted[id=r2]{Besides, the implementation of neural networks is relatively straightforward and requires little modification for different input dimensions. To that effect, }We have developed a suite of software libraries that make this method more easily accessible to other researchers without deep technical expertise in automatic differentiation or optimization. The code is accessible through
\begin{center}
\url{https://github.com/kailaix/ADCME.jl}
\end{center}

The remainder of this paper is organized as follows. We first introduce the problem setup in \Cref{sect:ps}, including the governing equations and the numerical scheme. Then in \Cref{sect:ml}, we present our general framework for combining FEM and neural networks, specifically its training process and the predictive process. After that, we briefly discuss its applicability and provide an accuracy analysis based on a model problem in \Cref{sect:aa}. In \Cref{sect:uq}, we present an approach for quantifying the predictive errors due to the heterogeneity of the material and the model approximations. Finally, we apply the framework to a simple linear equation, a multi-scale fiber-reinforced thin plate problem, and a highly nonlinear rubber membrane problem in \Cref{sect:app}. We conclude and discuss a possible generalization of the framework in \Cref{sect:conc}.

\section{Problem Setup}\label{sect:ps}

\subsection{Specific examples}

Let's consider a simple example to illustrate our problem. Consider a nonlinear Poisson equation with $\bx \in \RR^d$
\begin{equation}
	-\nabla \cdot \Big( \big( u(\bx) + \kappa(\bx) \big) \, \nabla u(\bx) \Big) = 0, \quad \bx \in \Omega \subset \RR^d
\end{equation}
with appropriate boundary conditions. 

We may formulate different learning problems. In the most difficult case, we may be learning a function $\Mcal_{\bt}$ (using a deep neural network with parameters $\bt$) that approximates:
\[
	\Mcal: (u, \bx) \mapsto 
	\big( u(\bx) + \kappa(\bx) \big) \, \nabla u(\bx)
\]
The input $u$ is a function in this case. In practice, the neural network $\Mcal_{\bt}$ takes as input a vector of discrete samplings $\{u(\bx_i)\}_{i=1,\ldots,n}$ (or some equivalent discrete representation), and $\bx$. The output is in $\RR^d$.

Other learning problems may involve
\[
	\Mcal: \bx \mapsto \kappa(\bx)
\]
or
\[
	\Mcal: \varepsilon(u)(\bx) \mapsto 
	\sigma(\bx)
\]
The latter case is shown in \autoref{sect:plate} for a 2D problem. The input $\varepsilon(u)(\bx)$ is a (3-dimensional) function of $\nabla u$, evaluated at location $\bx \in \RR^2$, and it represents the strain. The output $\sigma(\bx) \in \RR^3$ is the Cauchy-stress.

\subsection{Model Problems}
Consider a physical system described by static or steady partial differential equations
\begin{equation}
    \Pcal(u(\bx), \Mcal(u(\bx), \bx)) = \Fcal(u(\bx), \bx, p),\qquad  \bx\in \Omega
    \label{EQ:PDE}
\end{equation}
where the boundary conditions are excluded for brevity. The physical system is characterized by the generalized differential operator $\Pcal$ that defines a conservation relation or other type of balance law, the state variable $u(\bx)$ is the solution of the physical system on the space domain $\Omega$. The generalized differential operator $\Mcal$ defines the coarse-grained model, like constitutive relations in structure mechanics. And $\Fcal$ represents the external force term or other source terms, which depends on the parameter $p$.

The conservation relation $\Pcal$ is regarded as a fundamental law of nature. However, the modeling term $\Mcal(u(\bx), \bx)$, which contains empirical assumptions and simplifications brings uncertainties and imperfectness to the mathematical description. In the proposed approach, the modeling term $\Mcal(u(\bx), \bx)$ in \cref{EQ:PDE} is replaced by a space homogenized neural network $\Mcal_{\bt}(u(\bx))$, which could be designed to embed as much physical information and a priori knowledge as possible. The $\bt \in \R^{m}$ denotes the trainable parameters of the neural network. It is worth mentioning that different neural networks could be designed and applied to different computational areas when the physical properties of the problem vary in different areas.

\subsection{Discretization}
The discretizations and solution strategies of the conservation law $\Pcal$  have been well established.
When an appropriate discretization is applied to \cref{EQ:PDE},  it becomes
 \begin{equation}
 \bP(\bu, \Mcal_{\bt}(\bu)) - \bF(\bu, \bx, p) = \bm0
 \label{EQ:DISC_PDE}
 \end{equation}
where $\bu \in \R^{n}$ is the discrete state vector corresponding to the spatial discretization of $u$,
$\bP$ is the spatial discretization of the differential operator $\Pcal$, and  $\bF$ is the discrete external force vector.
In the present work, we mainly focus on solid mechanics applications, hence the FEM is applied to discretize the system.

Both the training process and the predictive process are based on the discrete \cref{EQ:DISC_PDE}. Bringing a given neural network model and the observed data $\bu$ into \cref{EQ:DISC_PDE},  the norm of the residual force is an indicator of the neural network model. Hence, we can train the neural network model by minimizing the norm of the residual force (known as the \textit{loss function}). And in the predictive process, \cref{EQ:DISC_PDE} is solved by the Newton method, which guarantees that predicted results satisfy the conservation laws. This marks the difference between current work and other data-driven paradigms~\cite{kirchdoerfer2016data, kirchdoerfer2017data}, which impose the conservation laws through constraints.


\section{Data-driven Approach}\label{sect:ml}
In this section, data-driven techniques, applied to represent and extract the unknown modeling term $ \Mcal(\bu, \bx)$ in \cref{EQ:PDE} are introduced.
The data-driven model  $\Mcal_{\bt}(u(\bx))$ is one kind of phenomenological models,  which avoid unaffordable computational cost paid for models derived from first principles.
And data-driven models, combined only correct domain knowledge with sufficient data, can discern the underlying patterns and structure.
Hence, they have been shown to outperform the other empirical models in some previous studies~\cite{ling2016reynolds,ling2016machine}.

\subsection{Neural Networks}\label{sect:nn}
The neural network itself is not an algorithm, but a framework to represent a complex model relating data inputs and data outputs.
The framework is composed of several connected layers. Each layer takes in the output $\bx$ of the previous layer,
transforms the inputs through an activation function $f(\bx)$ and outputs the result to the next layer.
A nonlinear layer is defined as $\mathbf{f}(\bx) = \sigma(\bW \bx + \mathbf{b})$, here $ \bW $ is the weight matrix and $\mathbf{b}$ is the bias. $\sigma$ is called the activation function, such as the identity function, \texttt{tanh} and the sigmoid function $\sigma(x) = \frac{1}{1+e^{-x}}$. Multi-layer neural networks are compositions of many such functions, which can be conveniently written as
\begin{equation}
    \mathbf{g}(\bx) = \mathbf{f}_1\cdot \mathbf{f}_2 \cdots \mathbf{f}_L (\bx)
\end{equation}
\replaced[id=r2]{It was shown in \cite{lorentz1996constructive} that a one-layer feed-forward neural network with infinite width can approximate any continuous functions on a compact domain (known as the ``universal approximation theorem'').}{Such framework features the so called ``universal approximation'' property, which states that a one-layer feed-forward neural network with sufficient number of
neurons can approximate any continuous functions on a compact subset, under mild assumptions on the activation functions~\cite{lorentz1996constructive}.}  Particularly, we have explicit approximation error bounds for one and two layer neural networks if the sigmoid activation functions are used~\cite{mhaskar1996neural, xu2018calibrating, 2019arXiv190107758X}.
Moreover, neural networks suffer less the curse of dimensionality for high dimensional problems~\cite{weinan2017deep} compared with polynomial approximations and can avoid Gibbs phenomenon for discontinuous problems under certain conditions~\cite{llanas2008constructive}. \added[id=r2]{However, we need to be cautious when interpreting the approximation theory: the universal approximation theorem only states the existence of a good solution but does not provide an efficient way to find it.  In fact, the corresponding optimization problem is usually non-convex and many local minima may exist, making the optimization computationally challenging. Additionally, the theorem says nothing about the efficiency of the approximation, i.e., the convergence rate regarding the width of the neural network. Nevertheless, those concerns do not prevent us from exploring the expressive power of neural networks.}
\deleted{These provide us a strict mathematical justification of approximating unknown functions or models in the physical systems by using neural networks.}
\replaced[id=r2]{In this work, }{Besides,} a detailed comparison between the neural network model and its counterparts, including piecewise linear functions, \replaced[id=r2]{radial basis functions, and radial basis function networks}{, and radial basis functions}, is presented in \cref{sect:cplrs_rbf}, which demonstrates the superior regularization and generalization properties of the neural network model.

In the present study, a \replaced[id=r1]{fully-connected}{three-layer} neural network is applied for computational efficiency to approximate the unknown functions, such as the nonlinear constitutive relation in this paper. \added{We have chosen the fully-connected neural network architecture since it is commonly used and easy to implement. Additionally, we do not use very deep or wide neural networks so that they are easier to train.}
However, when the physical properties of the underlying model are available, we can design special architecture to enforce the physical constraints or accelerate the computation. Besides, if the unknown function is complicated, deeper neural networks are preferred since they have been demonstrated numerically to be more expressive.

\subsection{Training Process}\label{sect:tp}
Most of the neural network training processes~\cite{ghaboussi1991knowledge, ellis1995stress, shen2005finite, ling2016machine, furukawa1998implicit, wang2018multiscale} are direct constitutive relation fitting, which rely on the cleaned input  and output data of the modeling term $\Mcal(u(\bx), \bx)$, like strain vs stress data.
But for complex materials or phenomena, the measurement or computation of high fidelity comprehensive input-output data of $\Mcal(u(\bx), \bx)$ would be challenging.
In the present work, the neural network model is trained by an end-to-end approach, namely using data $u(\bx)$ and the associated load conditions, measured by full-field measurement techniques
such as digital image correlation and the grid method~\cite{surrel1994moire}
or generated by high-fidelity numerical simulations.
We assume the data set contains pairs of $(\bu_i,\bF_i)$, $i=1$,$\ldots$, $N$, or only the external force and boundary conditions, which are enough to assemble the external force term $\bF_i$.
For most engineering applications, the data are limited, which are obtained by either high-fidelity simulations or experiments.
Therefore, the present training process should be suitable and effective with a small data set.
Thanks to the enormous richness of constitutive information contained in these data,
the training process takes about $\mathcal{O}(10)$ data pairs in the present applications.

By substituting the unknown constitutive relation $\Mcal(u(\bx), \bx)$ in \cref{EQ:PDE} with the neural network approximation $\Mcal_{\bt}(\bu)$, which takes discretized displacement vector $\bu\in \RR^n$ or its associated strain field as the input, and outputs the stress field, we can formulate the loss function as
\begin{equation}
    L(\bt) = \sum_{i=1}^{N} \left(\bP(\bu_i, \Mcal_{\bt}(\bu_i)) - \bF_i \right)^2
    \label{eq:loss}
\end{equation}
\added{which depends on the full-field data, assuming we have sufficient resolution to construct the associated derivatives.} Then $L(\bt)$ is minimized to obtain an optimal parameter estimator $\hat \bt\in\Theta$.

The optimization of the loss function \cref{eq:loss} to determining the weights $\hat\bt$ can be done by gradient descent methods, specifically the Limited-memory BFGS~(L-BFGS-B) method \added{\cite{byrd1995limited}} in the present work. \added{We use the line search routine in \cite{more1994line}, which attempts to enforce the Wolfe conditions \cite{byrd1995limited} by a sequence of polynomial interpolations. Note BFGS is applicable in our case since the data sets are typically small, otherwise the more commonly used stochastic gradient descent for training neural networks should be adopted. The choice of BFGS optimizer is well-motivated for scientific and engineering applications due to limited data.} Modern frameworks such as \texttt{TensorFlow} we adopt in the paper provide us a way of computing the gradients
$\nabla_{\bt} L(\bt)$ using reverse-mode automatic differentiation~(AD).  It applies symbolic differentiation at the elementary operation level. In AD, all numerical computations are ultimately compositions of a finite set of elementary operations for which derivatives are known, and combining the derivatives of the constituent operations through the chain rule gives the derivation of the overall composition. AD has forward-mode and reverse-mode. A thorough investigation of their properties is beyond the scope of this paper. In a nutshell, researchers only need to focus on the forward simulation. The differentiation and optimization parts are taken care of by the software. The traditional solver and the neural network are combined to fulfill the end-to-end training process. An illustrative example is presented in the Appendix.

Moreover, each evaluation of the loss function \cref{eq:loss} and its gradient does not require solving the linear or nonlinear system, where an expensive Newton's solver may be required in the latter case. Hence it is efficient even when the optimization needs thousands of steps to converge, which is the general case for large scale neural network optimization.

\subsection{Prediction Process}\label{sect:pp}
\replaced[id=r2]{For the prediction}{The prediction process is straightforward}, Newton's method with the load stepping is applied to solve \cref{EQ:DISC_PDE}. 
Neural network models can efficiently deliver the prediction values and their derivatives with respect to any input variables.
Therefore, traditional solvers are applied with minor changes in the model term query.
And we can compute the Jacobian of $\bP$ in \cref{EQ:DISC_PDE} as
\begin{equation}
    D_{\bu} \bP(\bu, \Mcal_{\bt}(\bu)) = \nabla_{\bu} \bP(\bu, \Mcal_{\bt}(\bu)) + \nabla_{\Mcal_{\bt}} \bP(\bu, \Mcal_{\bt}(\bu))\nabla_{\bu} \Mcal_{\bt}(\bu)
\end{equation}
where $D_{\bu}$ designates the partial derivative with respect to the displacement field $\bu$. $\nabla_{\bu} \Mcal_{\bt}(\bu) $ is obtained via automatic differentiation. The Jacobian can be used for the Newton's solver.

\begin{figure}[htbp] 
\centering
\includegraphics[width=0.7\textwidth,keepaspectratio]{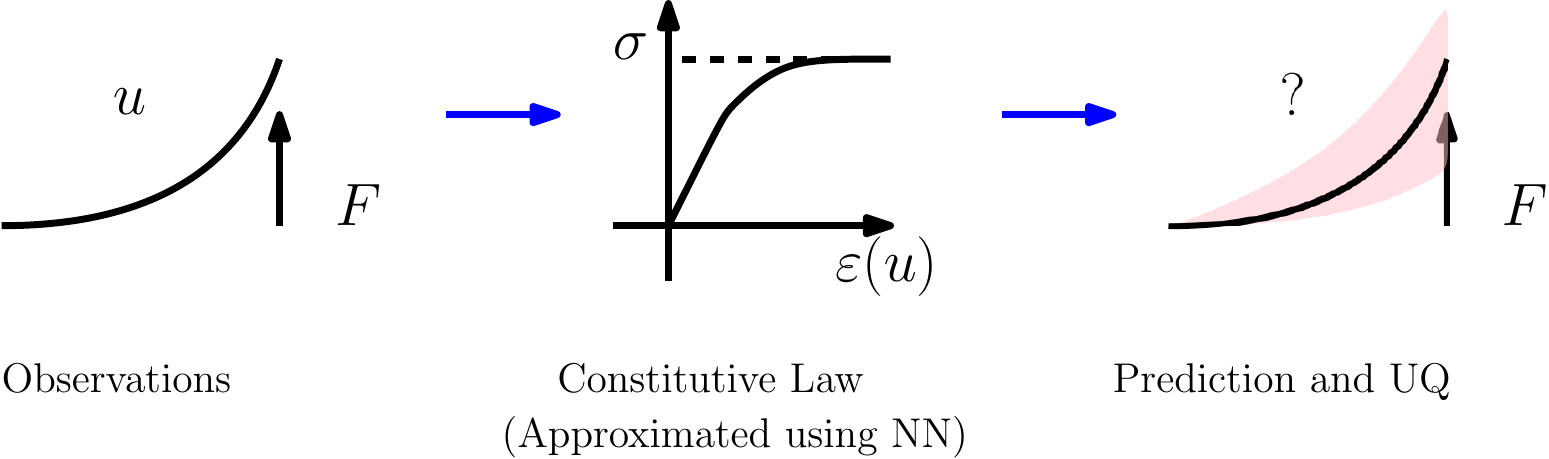}
\caption{Workflow for the predictive modeling and its uncertainty quantification.}
\label{fig:wf}
\end{figure}

The workflow described in \Cref{sect:nn,sect:tp,sect:pp}, from training the neural network model to predicting the material behaviors, is visualized in \cref{fig:wf}.

\section{Applicability and Accuracy Analysis}\label{sect:aa}
In this section, the conditions under which the aforementioned learning procedure is effective and the error bound of the predictive model learned by neural networks are discussed for a model problem.
Consider the 1D variable coefficient Poisson equation,
\begin{equation}
\begin{aligned}
&-\partial_x \big(\kappa(x) \partial_x u(x) \big)= f(x), \qquad 0 < x < 1 \\
&u(0) = u_0,\, u(1) = u_1
 \end{aligned}
 \label{EQ:MODEL_POISSON_PDE}
\end{equation}
Here $u_0$, $u_1$ are two numbers, $f(x)$ is the source function and the coefficient $\kappa(x)$ is approximated by a neural network $\kappa_{\bt}(x)$ parameterized by $\bt$.
In this case, the corresponding $\mathcal{M}_{\bt}(u,x)$ in the learning problem has the following physics format,
\begin{equation*}
\Mcal_{\bt}(u, x) =  \kappa_{\bt}(x)
\end{equation*}
We assume that $f(x)$ and $\kappa(x)$ have sufficient regularity so that $u(x)$ is also smooth.

{
The variational formulation of \cref{EQ:MODEL_POISSON_PDE} is discretized by the FEM on a uniform domain partition $\Tcal^h = \{0 = x_0 < x_1 < \cdots < x_{N_e} = 1\}$ with $h = x_j - x_{j-1}$, $j=1$, $2$, $\cdots$, $N_e$.
Let $\Ccal(P_1(\Tcal^h))$ denote the continuous piecewise linear function space on $\Tcal^h$,  a subspace of Sobolev space $H^{1}(0,1)$. The finite element formulation of \cref{EQ:MODEL_POISSON_PDE}
is given by: Find $u^h \in S^h = \{u | u\in \Ccal(P_1(\Tcal^h)),  u(0) = u_0,\, u(1) = u_1\}$  such that:
\begin{equation}
a(u^h , w^h) = \int_{0}^{1} \kappa(x) u^h_{,x} w^h_{,x} dx = \int_{0}^{1} f w^h dx  = (f, w^h)
\label{EQ: MODEL_POISSON_INT}
\end{equation}
holds for $\forall w^h \in V^h = \{w | w\in \Ccal(P_1(\Tcal^h)),  w(0) = 0,\, w(1) = 0\}$.  Applying one point Gaussian quadrature rule in each element,
the bilinear operator \cref{EQ: MODEL_POISSON_INT} leads to
\begin{equation}
\|a^h(u^h , w^h) - a(u^h, w^h)\| =  \| h\sum_{j = 1}^{N_e} \kappa(x_{j-1/2}) \frac{u^h(x_{j}) -  u^h(x_{j - 1}) }{h} \frac{w^h(x_{j}) -  w^h(x_{j - 1}) }{h} - a(u^h, w^h) \| \leq C_1 h^2
\label{EQ: MODEL_POISSON_INT_DISCRETE0}
\end{equation}
In the case that $w^h$ are local linear basis functions, the summation in \cref{EQ: MODEL_POISSON_INT_DISCRETE0} has at most two non-vanishing summands and the local error can be improved to $\mathcal{O}(h^3)$
\begin{equation}
    \| a^h(u^h , w^h)  - a(u^h, w^h) \| \leq C_1 h^3
\label{EQ: MODEL_POISSON_INT_DISCRETE}
\end{equation}

For the training process, we collect $N$ data pairs either from simulation or from experimental data,  $(u_1, f_1)$, $(u_2, f_2)$, $ ... $, $(u_N, f_N)$. The parameters $\bt$ for the neural network are updated by minimizing the loss function \cref{eq:loss}
\begin{equation}
L(\bt) = \sum_{i = 1}^N \| \bP_i -  \bF_i \|^2
\label{EQ: OBJECTIVE}
\end{equation}
here $\bP_i = \{a_{\bt}^h(u_i, \phi_1),  a_{\bt}^h(u_i, \phi_2), ..., a_{\bt}^h(u_i, \phi_{N_e - 1})\}$ and $\bF_i = \{(f_i, \phi_1),  (f_i, \phi_2), ..., (f_i, \phi_{N_e - 1})\}$ are assembled by
the FEM. And $\phi_i$ is the hat function at node $i$. $a_{\bt}^h$ is the discretized bilinear operator in \cref{EQ: MODEL_POISSON_INT_DISCRETE}  equipped with the neural network constitutive relation.
Assume that we are able to minimize the objective function so that the relative error for each term is bounded by $\mathcal{O}(\epsilon_0)$, i.e.,
\begin{equation}
 \frac{\| \bP_i -  \bF_i \|}{\|\bF_i\|} \leq C_2 \epsilon_0,  \quad  \forall i = 1,2,...,N
\label{EQ: OBJECTIVE_ERROR}
\end{equation}
Since $f$ is smooth, we have
\begin{equation}
    \|(\bF_i)_j \|= \|(f_i, \phi_j)\| \leq C_3 h,  \quad  \forall j  = 1,2,...N_e - 1 \textrm{ and } i = 1,2,...,N
\end{equation}
we have $\|\bF_i\|^2 \leq C_3^2 (N_e-1)h^2 \leq C_3^2 h$.
Therefore, on average, the optimization error of each component of \cref{EQ: OBJECTIVE_ERROR} satisfies
\begin{equation}
\|a_{\bt}^h(u_i , \phi_j)  -  (f_i, \phi_j)\| \approx \sqrt{\frac{\| \bP_i -  \bF_i \|^2}{N_e-1}} \leq  2C_2C_3 h\epsilon_0 ,  \quad  \forall j  = 1,2,...N_e - 1 \textrm{ and } i = 1,2,...,N
\label{EQ: ERROR_OPT}
\end{equation}

Plugging the data into \cref{EQ: MODEL_POISSON_INT}, and combining with  \cref{EQ: MODEL_POISSON_INT_DISCRETE} lead to
\begin{equation}
\| a^h(u_i , \phi_j)  - (f_i,  \phi_j) \| = \|a^h(u_i , \phi_j)  - a(u_i , \phi_j)\| \leq C_1h^3, \quad  \forall j  = 1,2,...,N_e - 1 \textrm{ and } i = 1,2,...,N
\label{EQ: ERROR_DISCRETE_LOCAL}
\end{equation}
Subtracting \cref{EQ: ERROR_OPT} from \cref{EQ: ERROR_DISCRETE_LOCAL}, we obtain
\begin{equation}
\|(a^h - a_{\bt}^h)(u_i , \phi_j)\|  \leq 2C_2C_3h\epsilon_0 + C_1h^3
\label{EQ: ERROR_NN}
\end{equation}
Bringing \cref{EQ: MODEL_POISSON_INT_DISCRETE} and the definition of the hat function into \eqref{EQ: ERROR_NN} leads to
\begin{equation}
\|(\kappa(x_{j - \frac{1}{2}}) -  \kappa_{\bt}(x_{j - \frac{1}{2}}))\frac{u_i(x_{j}) - u_i(x_{j - 1})}{h} -  (\kappa(x_{j + \frac{1}{2}}) -  \kappa_{\bt}(x_{j + \frac{1}{2}}))\frac{u_i(x_{j + 1}) - u_i(x_{j})}{h} \| \leq2C_2C_3h\epsilon_0 + C_1h^3
\label{EQ: ERROR_NN2}
\end{equation}
Consider another data pair $(u_k, f_k)$, we can obtain the same estimation as \cref{EQ: ERROR_NN2}
\begin{equation}
\|(\kappa(x_{j - \frac{1}{2}}) -  \kappa_{\bt}(x_{j - \frac{1}{2}}))\frac{u_k(x_{j}) - u_k(x_{j - 1})}{h} -  (\kappa(x_{j + \frac{1}{2}}) -  \kappa_{\bt}(x_{j + \frac{1}{2}}))\frac{u_k(x_{j + 1}) - u_k(x_{j})}{h}\| \leq 2C_2C_3h\epsilon_0 + C_1h^3
\label{EQ: ERROR_NN2_P}
\end{equation}
Combining \cref{EQ: ERROR_NN2}, \cref{EQ: ERROR_NN2_P} and matrix norm definition leads to
\begin{equation}
\Big\|
\begin{bmatrix}
    \kappa(x_{j - \frac{1}{2}}) -  \kappa_{\bt}(x_{j - \frac{1}{2}})\\
    \kappa(x_{j + \frac{1}{2}}) -  \kappa_{\bt}(x_{j + \frac{1}{2}})\\
\end{bmatrix}
\Big\|
\leq
\Big\|
\begin{bmatrix}
    \frac{u_i(x_{j}) - u_i(x_{j - 1})}{h}      & \frac{u_i(x_{j + 1}) - u_i(x_{j})}{h} \\
    \frac{u_k(x_{j}) - u_k(x_{j - 1})}{h}      & \frac{u_k(x_{j + 1}) - u_k(x_{j})}{h} \\
\end{bmatrix}^{-1}
\Big\|\cdot
\Big\|
\begin{bmatrix}
    2C_2C_3h\epsilon_0 + C_1h^3 \\
    2C_2C_3h\epsilon_0 + C_1h^3 \\
\end{bmatrix}
\Big\|
\label{EQ: ERROR_NN3}
\end{equation}
Through Taylor expansion of the data $u_i$ and $u_k$, we have
\begin{equation}
\det{\begin{vmatrix}
    \frac{u_i(x_{j}) - u_i(x_{j - 1})}{h}      & \frac{u_i(x_{j + 1}) - u_i(x_{j})}{h} \\
    \frac{u_k(x_{j}) - u_k(x_{j - 1})}{h}      & \frac{u_k(x_{j + 1}) - u_k(x_{j})}{h} \\
\end{vmatrix}}
 = h \left(u_i'(x_j)u_k''(x_j) - u_k'(x_j)u_i''(x_j)\right) + \mathcal{O}(h^2)
\end{equation}

The error bound of the neural network model at each Gaussian point solving by \cref{EQ: ERROR_NN3}  is given as
\begin{equation}
\begin{aligned}
&\|\kappa(x_{j - \frac{1}{2}}) -  \kappa_{\bt}(x_{j - \frac{1}{2}})\| \cdot | u_i'(x_j)u_k''(x_j) - u_k'(x_j)u_i''(x_j) |   \leq C_4(2C_2C_3\epsilon_0 + C_1h^2),\\
&\|\kappa(x_{j + \frac{1}{2}}) -  \kappa_{\bt}(x_{j + \frac{1}{2}})\| \cdot | u_i'(x_j)u_k''(x_j) - u_k'(x_j)u_i''(x_j) |   \leq C_4(2C_2C_3\epsilon_0 + C_1h^2)
\end{aligned}
\label{EQ: ERROR_BOUND}
\end{equation}
When both $\kappa$ and $\kappa_{\bt}$ are smooth enough,
through interpolation of \cref{EQ: ERROR_BOUND}, we can obtain the error bound on the interior point $x_{j - \frac{1}{2}} \leq x \leq x_{j + \frac{1}{2}}$ as follows
\begin{equation*}
\|\kappa(x) -  \kappa_{\bt}(x)\| \cdot | u_i'(x_j)u_k''(x_j) - u_k'(x_j)u_i''(x_j) |  \leq  C_4(2C_2C_3\epsilon_0 + C_1h^2) + C_5h^2
\end{equation*}
in other words,
\begin{equation}
\boxed{
    \|\kappa(x) -  \kappa_{\bt}(x)\|   \leq  \left(\frac{2C_2C_3C_4\epsilon_0 + (C_1C_4 + C_5)h^2}{\gamma_j}\right),\qquad \gamma_j = \sup_{i,k=1,2,\cdots, N} | u_i'(x_j)u_k''(x_j) - u_k'(x_j)u_i''(x_j) |
}
    \label{EQ: ERROR_BOUND2}
\end{equation}

The error bound in \cref{EQ: ERROR_BOUND2} reveals a quantitative relationship between optimization, discretization, and data. The error term consists of the optimization error $\mathcal{O}\left(\epsilon_0 \right)$  and the discretization error $\mathcal{O}\left(h^2\right)$, with constants $C_i$ independent of both $\epsilon_0$ and $h$. The errors are magnified by the reciprocal of the correlation of the data $\frac{1}{\gamma_j}$. If there are sufficient data, we can have a lower bound for the correlation term $\gamma_j$. In the limit case $h\rightarrow 0$ and $\varepsilon_0\rightarrow 0$, we obtain the convergence of $\kappa_{\bt}(x)$ to the true coefficient $\kappa(x)$. For optimization error dominant cases, the lesson is that increasing mesh resolution may not improve model learning, namely, fully-resolved meshes are not necessary for problems with coarse-grained models. Additionally, the error bound is an upper bound and is only valid under certain assumptions on the neural networks, optimization and model regularity.
}

In sum, the present model is applicable and effective, when the following conditions are simultaneously satisfied
\begin{itemize}
\item The neural network is consistent, namely with correct input features, the optimization error should tend to zero.
\item Both the underlying model $\kappa$ and the predicted model $\kappa_{\bt}$ are smooth enough to have bounded derivatives.
\item The data $u_{i}$ and $u_k$ should not be too correlated such that $| u_i'(x_j)u_k''(x_j) - u_k'(x_j)u_i''(x_j) |$ is small or vanishes. However, the issue can be resolved, when the data set is large enough so that there exist sufficient non-correlated observations.
\end{itemize}


\section{Uncertainty Quantification~(UQ)}\label{sect:uq}
Despite the error bound derived in \Cref{sect:aa},  the aforementioned neural network model contains uncertainties, due to the optimization error, incomplete input features of the neural network, data noise, and the homogenization error. Estimating the uncertainties of the FEM-neural network framework is critical for predictive modeling. There exists lots of prior work to quantify system uncertainties,
such as Monte Carlo methods, which rely on repeated random forward sampling,
polynomial chaos methods~\cite{xiu2002wiener,xiu2003modeling}, which determine the evolution of input uncertainty in a dynamical system through orthogonal polynomials,
and Bayesian procedures~\cite{rasmussen2004gaussian, kennedy2001bayesian},  which infer the posterior distribution of unknowns from existent data.
And more recently deep learning techniques such as Dropout~\cite{gal2016dropout} and DNN-based surrogate~\cite{TRIPATHY2018565} are applied to quantify the uncertainties in the neural networks.
In this section, we propose a UQ method specifically for quantifying the homogenization error, when the model $\Mcal_{\bt}(u)$ does not depend on $\bx$. This is similar to
the neural network error due to the incompleteness of its input features. Because the heterogeneity information, i.e.the coordinate $\bx$, is not incorporated in the neural network model in the present paper.

Solving the discretized governing equation~(\cref{EQ:DISC_PDE}), we have
\begin{equation}
        \hat \bu = \bu(\hat\bt, p)
        \label{EQ: UQ_PRED}
\end{equation}
where $p$ denotes the force load parameter.  The approximated constitutive relation model parameter $\hat\bt\in \RR^m$ is learned from data by minimizing \cref{eq:loss}. We have assumed a homogenized model, i.e., the constitutive relation is assumed to be space-invariant in the computational domain.
However, in reality, at each element or each Gaussian point, the constitutive relations are slightly different due to the heterogeneity of the material. Therefore, the true solution is
\begin{equation}
        \bu = \bu(\bt_1,\bt_2,\,\cdots\,\bt_{g}, p)
        \label{EQ: UQ_REAL}
\end{equation}
where $g$ denotes the total number of the Gaussian quadrature points over the whole computational domain, and $\bt_i\in \RR^m$ is the parameter associated to the constitutive model at the $i$-th Gaussian point.  The discrepancy between \cref{EQ: UQ_REAL} and \cref{EQ: UQ_PRED} is defined as the uncertainty derived from homogenization.

Taylor expansion of the difference between \cref{EQ: UQ_PRED} and \cref{EQ: UQ_REAL} at $\hat{\bt}$ is written as 
\begin{equation}
\resizebox{.9 \textwidth}{!} 
{$
        \Delta \bu = \hat{\bu} - \bu(\bt_1,\bt_2,\,\cdots\,\bt_{g}, p) = \bu(\hat \bt,\, \hat \bt,\,\cdots\, \hat \bt,\, p) - \bu(\bt_1,\bt_2,\,\cdots\,\bt_{g},\, p) \approx \frac{\partial \bu}{\partial (\bt_1,\bt_2,\,\cdots\,\bt_{g})} [\Delta \bt_1,\,\Delta \bt_2,\,\cdots\,\Delta \bt_{g} ]^{T}
$} \label{EQ: UQ_DIFF_F1}
\end{equation}
here $\Delta \bt_i  = \hat\bt - \bt_i\in \RR^m$, $i = 1$, $2$,$\cdots$, $g$, represent the constitutive relation model form uncertainties on each Gaussian point.
They are assumed to be independent and identical distributed~(i.i.d.) random variables.
And we assume the error has no bias, namely $\E\left[ \Delta \bt \right]= \mathbf{0}$. Its variance $\bSig_{\bt}$ is estimated from
$N$ training data by solving the following least square~(LSQ) problem,
\begin{equation}
\resizebox{.9 \textwidth}{!} 
{$
\begin{aligned}
     \min_{\bSig_{\bt}\succeq 0}  &\qquad \sum_{j=1}^N \sum_{i=1}^{n} \left( (\Delta u_j^i)^2 - \frac{\partial u_j^i}{\partial (\bt_1,\bt_2,\,\cdots\,\bt_{g})} \E [\Delta \bt_1,\,\Delta \bt_2,\,\cdots\,\Delta \bt_{g} ]^{T}  [\Delta \bt_1,\,\Delta \bt_2,\,\cdots\,\Delta \bt_{g} ]\Big( \frac{\partial u_j^i}{\partial (\bt_1,\bt_2,\,\cdots\,\bt_{g})} \Big)^T\right)^2 \\
 = & \sum_{j=1}^N  \sum_{i=1}^{n} \left( (\Delta u^i_j)^2 - \frac{\partial u_j^i}{\partial (\bt_1,\bt_2,\,\cdots\,\bt_{g})} \textrm{diag} \{ \bSig_{\bt},\,\bSig_{\bt},\,\cdots\,\bSig_{\bt}\}\Big( \frac{\partial u_j^i}{\partial (\bt_1,\,\bt_2,\,\cdots\,\bt_{g})} \Big)^T\right)^2
\label{EQ: UQ_DIFF_F2}
\end{aligned}
 $}
\end{equation}
 $\Delta u^i_j$ is the $i$-th component of the error vector   $\bu(\hat\bt, p_j) - \bu_j$, here $\bu_j$ is the $j$-th true solution~(or observation), and  $u_j^i$ is the $i$-th component of \cref{EQ: UQ_PRED} with force load $p_j$.

However, in most cases, the parameter number is large, the estimation of the variance matrix $\bSig_{\bt}\in \RR^{m\times m}$ needs a large amount of data.
Based on the idea of active subspace methods proposed in~\cite{constantine2014active},
the random variable $\Delta \bt$ is restricted to a low-dimensional subspace of $\R^{m}$, represented by an associated matrix denoted here by $\bW \in R^{m\times k}$, whose dimension $k$ is an order of magnitude smaller than $m$. The subspace is constructed by recovering an orthogonal projection matrix $\bW$ obtained through the orthogonalization of a linear manifold spanned by the gradients of quantities of interest~(QoIs)
\begin{equation}
    \mathrm{span}\{ \nabla_{\bt} J^1(\bu(\hat\bt,\, p)),\, \nabla_{\bt} J^2 (\bu(\hat\bt,\, p)),\, \cdots, \nabla_{\bt} J^k(\bu(\hat\bt,\, p)) \}
\end{equation}
here $J^{i}, i = 1,2,\,\cdots,\, m$ are QoIs.  The random difference at each Gaussian point for a given $p$ is modeled as
\begin{equation}
\Delta \bt_i  =  \blam_i  \bW^T \qquad i=1,2,\cdots, g
\label{EQ:MOD}
\end{equation}
It is worth mentioning $\bW$ is force load $p$ dependent, but we omit $p$ in the notation for brevity. The vector  of  reduced  coordinates  $\blam_i \in \R^{k}$ is a zero-mean i.i.d. random variable with a variance matrix $\bSig_{\blam}$.
Bringing \cref{EQ:MOD} into \cref{EQ: UQ_DIFF_F1} leads to
\begin{equation}
        \Delta u^i_j \approx \frac{\partial u^i_j}{\partial (\blam_1,\blam_2,\,\cdots\,\blam_{g})}
[\blam_1,\blam_2,\,\cdots\,\blam_{g}]^{T}
        \label{EQ: UQ_DIFF_R1}
\end{equation}
here
\begin{equation}
       \frac{\partial u^i_j}{\partial (\blam_1,\blam_2,\,\cdots\,\blam_{g})}  = \frac{\partial u^i_j}{\partial (\bt_1,\bt_2,\,\cdots\,\bt_{g})}
\textrm{diag}\{\bW,\,\bW,\,\cdots,\,\bW\}
\end{equation}
The variance matrix $\bSig_{\blam}\in \RR^{k\times k}$ is approximated by solving the following least square problem
\begin{equation}
\begin{aligned}
       \min_{\bSig_{\blam}\succeq 0}  \qquad  \sum_{j=1}^N \sum_{i=1}^n \left( (\Delta u^i_j)^2 - \frac{\partial u^i_j}{\partial (\blam_1,\blam_2,\,\cdots\,\blam_{g})}
\textrm{diag}\{\bSig_{\blam},\,\bSig_{\blam},\,\cdots,\,\bSig_{\blam}\}  \Big(\frac{\partial u^i_j}{\partial (\blam_1,\blam_2,\,\cdots\,\blam_{g})} \Big)^T\right)^2
\label{EQ: UQ_DIFF_T2}
\end{aligned}
\end{equation}

In the present framework, the low-dimensional subspace is chosen to be one-dimensional~($k$ = 1) and we use $\lambda\in\RR$ to denote the reduced coordinate. The only QoI, $J$, is taken to be the maximum principal stress. The parameter is approximated in the one-dimensional subspace,
\begin{equation}
\Delta \bt = \lambda \frac{\nabla_{\bt} J}{ \| \nabla_{\bt} J \|}
\end{equation}
The normalization is necessary since $J$ can be quite different in scales for different external loads.

\Cref{EQ: UQ_DIFF_T2} can be further simplified as
\begin{equation}
\begin{aligned}
         \min_{\bSig_{\blam}\succeq 0}  \qquad \sum_{j=1}^N \sum_{i=1}^n \left((\Delta u^i_j)^2 - \bSig_{\blam} c_j^i\right)^2
\label{EQ: UQ_DIFF_R3}
\end{aligned}
\end{equation}
here $c_j^i = \frac{\partial u^i_j}{\partial (\blam_1,\blam_2,\,\cdots\,\blam_{g})}  \Big(\frac{\partial u^i_j}{\partial (\blam_1,\blam_2,\,\cdots\,\blam_{g})} \Big)^T$. And normalizing the summand in \cref{EQ: UQ_DIFF_R3} by $\frac{1}{(c_j^i)^2}$, when $c_j^i > 0$, can improve the least square estimation.

Based on the Chebyshev's inequality, the LSQ estimated variance $\hat \bSig_{\blam}$ satisfies
\begin{equation}
\boxed{
P(|\bSig_{\lambda} - \hat \bSig_{\lambda}| \geq \epsilon) \leq \frac{\mathrm{Var}(\bSig_{\lambda})}{\epsilon^2 Nn}
}
\end{equation}
Therefore, when $Nn$ is large enough, the estimation $\hat \bSig_{\blam}$ obtained by the LSQ converges to $\bSig_{\blam}$ with high probability.

During the prediction process, the Monte Carlo sampling method is applied to compute the confidence interval
 \begin{equation}
     \bu_{\mathrm{pred}}(p)  \approx \bu(\hat \bt, p) + \frac{\partial \bu}{\partial (\blam_1\
,\blam_2,\,\cdots\,\blam_{g})} [\blam_1,\blam_2,\,\cdots\,\blam_{g}]^{T}
\label{EQ:MC}
 \end{equation}
where $\blam_i$  is generated as a Gaussian random variable, $\mathcal{N}(0, \hat \bSig_{\lambda})$. It is worth mentioning that the sampling process does not require repeated solving the
forward problem, which is efficient for large sampling.

\begin{remark}  { When the measurements are noisy, the $j$-th observed data is actually $\bu^{ob}_j = \bu_j + \bm{e}_j$ and therefore
$$\underbrace{\bu^{\mathrm{pred}}_j-\bu_j}_{\mbox{model error}} = \underbrace{\bu^{\mathrm{pred}}_j-\bu^{ob}_j}_{\mbox{total error}} + \underbrace{\bm{e}_j}_{\mbox{measurement error}} $$
 Here $\bm{e}_j$ denotes the measurement error and $\bu^{\mathrm{pred}}_j$ is the prediction.  Assuming that the measurement error and total error are uncorrelated and unbiased, we have
 $$
 \mathbb{E}(||\bu^{\mathrm{pred}}_j-\bu_j||^2) = \mathbb{E}||\bu^{\mathrm{pred}}_j-\bu^{ob}_j||^2 + \mathbb{E}||\bm{e}_j||^2 
 $$
 The LSQ problem~\cref{EQ: UQ_DIFF_F2} needs to be modified and may be rewritten as
\begin{equation}
\label{EQ:UQ_MEASURE}
\resizebox{.9 \textwidth}{!} 
{$
\begin{aligned}
     \min_{\bSig_{\bt}\succeq 0} \qquad \sum_{j=1}^N  \sum_{i=1}^{n} \left( \underbrace{(\Delta u^{ob,i}_j)^2}_{\mbox{total error variance}} + \underbrace{\E(e_j^i)^2}_{\mbox{measurement error variance}}- \underbrace{\frac{\partial u_j^i}{\partial (\bt_1,\bt_2,\,\cdots\,\bt_{g})} \textrm{diag} \{ \bSig_{\bt},\,\bSig_{\bt},\,\cdots\,\bSig_{\bt}\}\Big( \frac{\partial u_j^i}{\partial (\bt_1,\,\bt_2,\,\cdots\,\bt_{g})}}_{\mbox{model error variance}} \Big)^T 
     \right)^2
\end{aligned}
 $}
\end{equation}
$\Delta u^{ob,i}_j$ and $e_j^i$ are the $i$-th component of the observed state vector and the error vector. Neglecting the measurement errors in \cref{EQ:UQ_MEASURE} may lead to under-estimating the model form uncertainty $\bSig_{\bt}$.
After that, the following UQ procedure remains the same. For this approach to be applicable, $\E(e_j^i)^2$ must be provided. 
}
\end{remark}

\section{Applications}\label{sect:app}

In this section, we {first present a simple linear equation case for verifying the analysis in~\cref{sect:aa} and exploring neural network architectures.} Then we present numerical results from solid mechanics for the proposed ``small-data''-driven predictive modeling procedure: a multi-scale fiber-reinforced plate problem and a highly nonlinear rubbery membrane problem.
\begin{itemize}
    \item The first problem serves as a proof-of-concept example for the end-to-end approach, where we calibrate the linear fourth-order stiffness tensor in the constitutive relation. It can also be viewed as a demonstration of guess-then-fit approach: we first guess that the material is subject to linear constitutive relation and then we fit the parameters from data.
    \item The second problem tackles the nonlinear constitutive relation with neural networks. And the strength of the neural network approach is explored in a thorough comparison with the piecewise linear functions~(PL), \replaced[id=r2]{radial basis functions~(RBF), and radial basis function networks (RBFN)}{and radial basis functions~(RBF)}. We show that in this case, PL and RBF are either overfitting with large degrees of freedoms or under-fitting with small degrees of freedoms and are susceptible to noise. Meanwhile, neural networks generalize well and are quite robust to noise.
\end{itemize}
In both problems, the training data and test data of the displacement field are generated numerically, the underlying constitutive models are chosen to be space-invariant or space-varying, i.e., containing random noise. The predicted results and the corresponding confidence interval on the test data are reported.

\subsection{Convergence Analysis}\label{sect:ca}

{As a numerical demonstration of our theoretical analysis \mbox{\Cref{EQ: ERROR_BOUND2}}, we consider a simple linear equation so the only error sources come from discretization and optimization (for nonlinear equations we consider in the next few examples, another source of errors include the Newton-Raphson solver; in those cases, it is harder to analyze the sources of error and convergence rates).
\begin{equation}\label{equ:numerical_model}
  \kappa (x) \Delta u(x) = f(x),\ x\in (0,1) \quad u(0)=u(1)=0
\end{equation}
The exact $u(x)$, $\kappa(x)$ and $f(x)$ are given by 
\[
	u(x) = x(1-x)\quad\kappa(x) = \frac{1}{1+x^2}\quad f(x)= -\frac{2}{1+x^2}
\]
We discretize the interval $[0,1]$ by $0=x_0<x_1<\ldots<x_{N_e}=1$, $x_{i+1}-x_i=h=\frac{1}{N_e}$. We assume that the data $\{u(x_i)\}_{i=0}^{N_e}$ are given and $\kappa(x)$ is unknown. $\kappa(x)$ is approximated by a fully-connected neural network $\kappa_{\bt}(x)$ with 3 hidden layers with 20 neurons for each layer, where $\bt$ is the weights and biases. The activation function is \texttt{tanh}. The last layer is a linear layer with a scalar output. The loss function and the corresponding $\mathbf{P}_i$ and $\mathbf{F}_i$ are given by 
\[
L(\theta) = \sum_{i=1}^{N_e-1} \|\mathbf{P}_i-\mathbf{F}_i\|^2,\quad \mathbf{P}_i = \kappa_{\bt}(x_i) \frac{u_{i+1}-2u_i+u_{i-1}}{h^2} ,\quad \mathbf{F}_i = f(x_i)\quad i = 1,2,\ldots, N_e - 1
\]
We consider $N_e=3,4,5,8,10,13,15,18,21$ and use the \texttt{L-BFGS-B} optimizer with the tolerance $10^{-12}$ and the maximum iteration $15000$ to reduce $\varepsilon_0$ as much as possible. Additionally, we generate $N_t = 1000$ points $\{x_i\}_{i=1}^{N_t}$ uniformly from $[0,1]$ as the test set (the test set is fixed for evaluating all computed $\bt$) and use
\[
\mathrm{E}_{\bt} =\sqrt{\frac{\sum_{i=1}^{N_t} |\kappa(x_i)-\kappa_{\bt}(x_i)|_2^2}{N_t}}
\]
as an estimation for $\|\kappa(x)-\kappa_{\bt}(x)\|_2$. 

The estimation for $\bt$ depends on many factors such as neural network architectures, initialization and optimization algorithms. To access the effect of discretization, we fix the neural network architecture and train the neural network with 10 different sets of initialized weights and biases. \mbox{\Cref{fig:boxplot}} shows the box plots of $E_{\bt}$ for different $N_e$ and initializations.  We see that the error keeps decreasing up to $N_e=8$. At this stage, the discretization error dominates. Note our analysis in \mbox{\Cref{EQ: ERROR_BOUND}} only gives an upper bound and the observed convergence rate is better than the predicted rate. When $N_e\geq 8$, the error $E_{\bt}$ flattens  because of the optimization error $\mathcal{O}(\varepsilon_0/\gamma_j)$. And $\varepsilon_0$ depends not only on the convergence of the optimizer but also the expressive power of the neural network. 

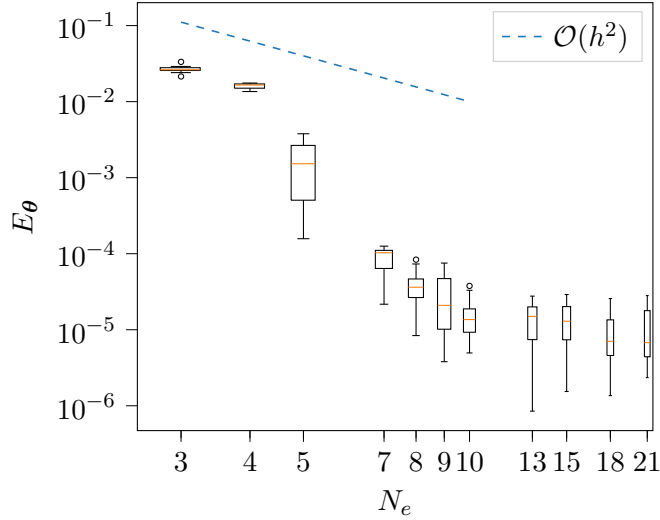
\begin{figure}[htbp] 
\centering
\scalebox{1.0}{
\begin{tikzpicture}

\definecolor{color0}{rgb}{1,0.498039215686275,0.0549019607843137}
\definecolor{color1}{rgb}{0.12156862745098,0.466666666666667,0.705882352941177}

\begin{axis}[
legend cell align={left},
legend style={draw=white!80.0!black},
log basis x={10},
log basis y={10},
minor xtick={3,4,5,7,8,9,10,13,15,18,21},
tick align=outside,
tick pos=left,
x grid style={white!69.01960784313725!black},
xlabel={$N_e$},
xmin=2.5, xmax=21.5,
xmode=log,
xtick style={color=black},
xtick={3,4,5,7,8,9,10,13,15,18,21},
xticklabels={3,4,5,7,8,9,10,13,15,18,21},
y grid style={white!69.01960784313725!black},
ylabel={\(E_{\bt}\)},
ymin=4.71524584441604e-07, ymax=0.200253017576893,
ymode=log,
ytick style={color=black},
ytick={1e-08,1e-07,1e-06,1e-05,0.0001,0.001,0.01,0.1,1,10},
yticklabels={\(\displaystyle {10^{-8}}\),\(\displaystyle {10^{-7}}\),\(\displaystyle {10^{-6}}\),\(\displaystyle {10^{-5}}\),\(\displaystyle {10^{-4}}\),\(\displaystyle {10^{-3}}\),\(\displaystyle {10^{-2}}\),\(\displaystyle {10^{-1}}\),\(\displaystyle {10^{0}}\),\(\displaystyle {10^{1}}\)}
]
\addplot [black, forget plot]
table {%
2.75 0.0256646733045546
3.25 0.0256646733045546
3.25 0.0279331879374074
2.75 0.0279331879374074
2.75 0.0256646733045546
};
\addplot [black, forget plot]
table {%
3 0.0256646733045546
3 0.0240906344186921
};
\addplot [black, forget plot]
table {%
3 0.0279331879374074
3 0.0289972896209899
};
\addplot [black, forget plot]
table {%
2.875 0.0240906344186921
3.125 0.0240906344186921
};
\addplot [black, forget plot]
table {%
2.875 0.0289972896209899
3.125 0.0289972896209899
};
\addplot [black, mark=*, mark size=1, mark options={solid,fill opacity=0}, only marks, forget plot]
table {%
3 0.0214202260427101
3 0.0333447728845994
};
\addplot [black, forget plot]
table {%
3.75 0.0150200991820031
4.25 0.0150200991820031
4.25 0.017089494630831
3.75 0.017089494630831
3.75 0.0150200991820031
};
\addplot [black, forget plot]
table {%
4 0.0150200991820031
4 0.0135704628453693
};
\addplot [black, forget plot]
table {%
4 0.017089494630831
4 0.0175624926369086
};
\addplot [black, forget plot]
table {%
3.875 0.0135704628453693
4.125 0.0135704628453693
};
\addplot [black, forget plot]
table {%
3.875 0.0175624926369086
4.125 0.0175624926369086
};
\addplot [black, forget plot]
table {%
4.75 0.000505301009519956
5.25 0.000505301009519956
5.25 0.00265534345892885
4.75 0.00265534345892885
4.75 0.000505301009519956
};
\addplot [black, forget plot]
table {%
5 0.000505301009519956
5 0.000157502743240489
};
\addplot [black, forget plot]
table {%
5 0.00265534345892885
5 0.00378084785859267
};
\addplot [black, forget plot]
table {%
4.875 0.000157502743240489
5.125 0.000157502743240489
};
\addplot [black, forget plot]
table {%
4.875 0.00378084785859267
5.125 0.00378084785859267
};
\addplot [black, forget plot]
table {%
6.75 6.39854205700487e-05
7.25 6.39854205700487e-05
7.25 0.000110651334477929
6.75 0.000110651334477929
6.75 6.39854205700487e-05
};
\addplot [black, forget plot]
table {%
7 6.39854205700487e-05
7 2.16629833152249e-05
};
\addplot [black, forget plot]
table {%
7 0.000110651334477929
7 0.000125660321089949
};
\addplot [black, forget plot]
table {%
6.875 2.16629833152249e-05
7.125 2.16629833152249e-05
};
\addplot [black, forget plot]
table {%
6.875 0.000125660321089949
7.125 0.000125660321089949
};
\addplot [black, forget plot]
table {%
7.75 2.6506480296019e-05
8.25 2.6506480296019e-05
8.25 4.65366703908065e-05
7.75 4.65366703908065e-05
7.75 2.6506480296019e-05
};
\addplot [black, forget plot]
table {%
8 2.6506480296019e-05
8 8.37958018127692e-06
};
\addplot [black, forget plot]
table {%
8 4.65366703908065e-05
8 7.3277842605765e-05
};
\addplot [black, forget plot]
table {%
7.875 8.37958018127692e-06
8.125 8.37958018127692e-06
};
\addplot [black, forget plot]
table {%
7.875 7.3277842605765e-05
8.125 7.3277842605765e-05
};
\addplot [black, mark=*, mark size=1, mark options={solid,fill opacity=0}, only marks, forget plot]
table {%
8 8.35747207516041e-05
};
\addplot [black, forget plot]
table {%
8.75 1.01594376624898e-05
9.25 1.01594376624898e-05
9.25 4.70555701419626e-05
8.75 4.70555701419626e-05
8.75 1.01594376624898e-05
};
\addplot [black, forget plot]
table {%
9 1.01594376624898e-05
9 3.80146799703664e-06
};
\addplot [black, forget plot]
table {%
9 4.70555701419626e-05
9 7.53252507497593e-05
};
\addplot [black, forget plot]
table {%
8.875 3.80146799703664e-06
9.125 3.80146799703664e-06
};
\addplot [black, forget plot]
table {%
8.875 7.53252507497593e-05
9.125 7.53252507497593e-05
};
\addplot [black, forget plot]
table {%
9.75 9.2671946439959e-06
10.25 9.2671946439959e-06
10.25 1.88182825800046e-05
9.75 1.88182825800046e-05
9.75 9.2671946439959e-06
};
\addplot [black, forget plot]
table {%
10 9.2671946439959e-06
10 4.9520773547541e-06
};
\addplot [black, forget plot]
table {%
10 1.88182825800046e-05
10 3.30338177781013e-05
};
\addplot [black, forget plot]
table {%
9.875 4.9520773547541e-06
10.125 4.9520773547541e-06
};
\addplot [black, forget plot]
table {%
9.875 3.30338177781013e-05
10.125 3.30338177781013e-05
};
\addplot [black, mark=*, mark size=1, mark options={solid,fill opacity=0}, only marks, forget plot]
table {%
10 3.77128571227568e-05
};
\addplot [black, forget plot]
table {%
12.75 7.42023875426317e-06
13.25 7.42023875426317e-06
13.25 1.99733072304513e-05
12.75 1.99733072304513e-05
12.75 7.42023875426317e-06
};
\addplot [black, forget plot]
table {%
13 7.42023875426317e-06
13 8.49817988065095e-07
};
\addplot [black, forget plot]
table {%
13 1.99733072304513e-05
13 2.76709468278741e-05
};
\addplot [black, forget plot]
table {%
12.875 8.49817988065095e-07
13.125 8.49817988065095e-07
};
\addplot [black, forget plot]
table {%
12.875 2.76709468278741e-05
13.125 2.76709468278741e-05
};
\addplot [black, forget plot]
table {%
14.75 7.37655563147434e-06
15.25 7.37655563147434e-06
15.25 2.01635804382856e-05
14.75 2.01635804382856e-05
14.75 7.37655563147434e-06
};
\addplot [black, forget plot]
table {%
15 7.37655563147434e-06
15 1.54095124736246e-06
};
\addplot [black, forget plot]
table {%
15 2.01635804382856e-05
15 2.90467106540804e-05
};
\addplot [black, forget plot]
table {%
14.875 1.54095124736246e-06
15.125 1.54095124736246e-06
};
\addplot [black, forget plot]
table {%
14.875 2.90467106540804e-05
15.125 2.90467106540804e-05
};
\addplot [black, forget plot]
table {%
17.75 4.5861483224199e-06
18.25 4.5861483224199e-06
18.25 1.34899560304506e-05
17.75 1.34899560304506e-05
17.75 4.5861483224199e-06
};
\addplot [black, forget plot]
table {%
18 4.5861483224199e-06
18 1.35975644398458e-06
};
\addplot [black, forget plot]
table {%
18 1.34899560304506e-05
18 2.57066062927908e-05
};
\addplot [black, forget plot]
table {%
17.875 1.35975644398458e-06
18.125 1.35975644398458e-06
};
\addplot [black, forget plot]
table {%
17.875 2.57066062927908e-05
18.125 2.57066062927908e-05
};
\addplot [black, forget plot]
table {%
20.75 4.40986763896311e-06
21.25 4.40986763896311e-06
21.25 1.78257106196534e-05
20.75 1.78257106196534e-05
20.75 4.40986763896311e-06
};
\addplot [black, forget plot]
table {%
21 4.40986763896311e-06
21 2.34646277535136e-06
};
\addplot [black, forget plot]
table {%
21 1.78257106196534e-05
21 2.81225380002818e-05
};
\addplot [black, forget plot]
table {%
20.875 2.34646277535136e-06
21.125 2.34646277535136e-06
};
\addplot [black, forget plot]
table {%
20.875 2.81225380002818e-05
21.125 2.81225380002818e-05
};
\addplot [semithick, color1, dashed]
table {%
3 0.111111111111111
5 0.04
8 0.015625
10 0.01
};
\addlegendentry{$\mathcal{O}(h^{2})$}
\addplot [color0, forget plot]
table {%
2.75 0.0262364131332296
3.25 0.0262364131332296
};
\addplot [color0, forget plot]
table {%
3.75 0.0165198290369848
4.25 0.0165198290369848
};
\addplot [color0, forget plot]
table {%
4.75 0.00152452581791111
5.25 0.00152452581791111
};
\addplot [color0, forget plot]
table {%
6.75 0.000102820901879767
7.25 0.000102820901879767
};
\addplot [color0, forget plot]
table {%
7.75 3.62793700852124e-05
8.25 3.62793700852124e-05
};
\addplot [color0, forget plot]
table {%
8.75 2.09250871011599e-05
9.25 2.09250871011599e-05
};
\addplot [color0, forget plot]
table {%
9.75 1.3612508205609e-05
10.25 1.3612508205609e-05
};
\addplot [color0, forget plot]
table {%
12.75 1.49822800824154e-05
13.25 1.49822800824154e-05
};
\addplot [color0, forget plot]
table {%
14.75 1.29442105819132e-05
15.25 1.29442105819132e-05
};
\addplot [color0, forget plot]
table {%
17.75 7.05822777945854e-06
18.25 7.05822777945854e-06
};
\addplot [color0, forget plot]
table {%
20.75 6.76373079898673e-06
21.25 6.76373079898673e-06
};
\end{axis}

\end{tikzpicture}}
\caption{Convergence plot for $E_{\bt}$. The error decreases until $N_e \sim8$. This indicates that the dominant error switches from the discretization error to the optimization/neural network approximation error at that point. Note our analysis only gives an upper bound and the observed convergence rate is better than the predicted rate $\mathcal{O}(h^2)$.}
\label{fig:boxplot}
\end{figure}

We also compare the effects of neural network width, depth, and activation functions. \mbox{\Cref{tab:nn}} shows the test errors for different configurations. For all the test cases, we let $N_e=50$ and run the optimization with 10 different sets of random initialization weights and biases for neural networks and report the mean and standard deviation of the errors. We can see that \texttt{tanh} is a better choice of activations among \texttt{tanh}, ReLU~(rectified linear unit)~\cite{nair2010rectified}, ELU~(exponential linear unit)~\cite{clevert2015fast}, SELU~(scaled exponential linear unit)~\cite{klambauer2017self}, and leaky ReLU~(leaky rectified linear unit)~\cite{he2015delving}. Given the activation function \texttt{tanh}, deeper neural networks achieve better accuracy compared to shallow ones for fixed width. Increasing width alone does not necessarily improve the accuracy, which we attribute to the optimization error.

\begin{table}[hbtp]
\begin{tabular}{c|c|ccc}
\toprule
& Width $\backslash$ Depth &1& 3 & 8 \\
\midrule
\multirow{4}{*}{\texttt{tanh}}   &2& 0.001553${}_{\pm0.0005344}$&$1.004\times 10^{-5}{}_{\pm8.487\times 10^{-6}}$&$4.179\times 10^{-6}{}_{\pm2.868\times 10^{-6}}$\\
&10& 0.03256${}_{\pm0.06785}$&$1.522\times 10^{-5}{}_{\pm1.58\times 10^{-5}}$&$1.232\times 10^{-5}{}_{\pm8.536\times 10^{-6}}$\\
&20& 0.04895${}_{\pm0.07753}$&$4.47\times 10^{-5}{}_{\pm6.861\times 10^{-5}}$&$8.308\times 10^{-6}{}_{\pm6.534\times 10^{-6}}$\\
&40& 0.1452${}_{\pm0.05092}$&0.09476${}_{\pm0.08132}$&$3.829\times 10^{-5}{}_{\pm6.723\times 10^{-5}}$\\
\midrule
\multirow{4}{*}{ReLU}    &2& 0.3059${}_{\pm0.265}$&0.4206${}_{\pm0.3445}$&0.3886${}_{\pm0.2726}$\\
&10& 0.1583${}_{\pm0.005901}$&0.2677${}_{\pm0.1293}$&0.6164${}_{\pm0.3752}$\\
&20& 0.1577${}_{\pm0.00091}$&0.1727${}_{\pm0.02352}$&0.5958${}_{\pm0.4045}$\\
&40& 0.1581${}_{\pm1.902\times 10^{-5}}$&0.1646${}_{\pm0.001314}$&0.3616${}_{\pm0.2365}$\\
\midrule
\multirow{4}{*}{ELU}    &2& 19.08${}_{\pm59.71}$&1.004${}_{\pm1.595}$&0.3506${}_{\pm0.3434}$\\
&10& 1.303${}_{\pm1.537}$&0.3694${}_{\pm0.4182}$&0.2974${}_{\pm0.2928}$\\
&20& 0.9138${}_{\pm1.025}$&0.351${}_{\pm0.2456}$&0.24${}_{\pm0.1018}$\\
&40& 0.5231${}_{\pm0.6619}$&0.3315${}_{\pm0.2392}$&0.9213${}_{\pm0.7828}$\\
\midrule
\multirow{4}{*}{SELU}    &2& 3.495${}_{\pm5.669}$&0.182${}_{\pm0.05718}$&0.2356${}_{\pm0.04606}$\\
&10& 1.153${}_{\pm1.099}$&0.4626${}_{\pm0.3673}$&0.5026${}_{\pm0.1774}$\\
&20& 0.6917${}_{\pm0.6361}$&0.5693${}_{\pm0.3387}$&0.6294${}_{\pm0.2241}$\\
&40& 0.299${}_{\pm0.1142}$&0.5105${}_{\pm0.2815}$&0.6738${}_{\pm0.1352}$\\
\midrule
\multirow{4}{*}{leaky ReLU}    &2& 0.6917${}_{\pm0.8444}$&0.3127${}_{\pm0.1943}$&0.3601${}_{\pm0.3482}$\\
&10& 4.006${}_{\pm8.434}$&0.6907${}_{\pm0.4909}$&0.7305${}_{\pm0.4462}$\\
&20& 0.2423${}_{\pm0.1325}$&1.234${}_{\pm1.026}$&0.7996${}_{\pm0.6328}$\\
&40& 0.3648${}_{\pm0.2608}$&0.2935${}_{\pm0.2009}$&0.4595${}_{\pm0.3247}$\\
\bottomrule
\end{tabular}
\caption{Error mean and standard deviation for different neural network architectures, including \texttt{tanh}, ReLU~(rectified linear unit), ELU~(exponential linear unit), SELU~(scaled exponential linear unit), and leaky ReLU~(leaky rectified linear unit). The coefficient of leakage in the leaky ReLU is 0.2. The smaller $\pm$ numbers denote the standard deviation. The width is the number of activation nodes in each hidden layer. The depth is the number of hidden layers.}
\label{tab:nn}
\end{table}

}

\subsection{Fiber Reinforced Plate}
\label{sect:plate}
\begin{figure}[htbp]
\centering
\begin{tikzpicture}
  \draw (0,0) -- (7,0) -- (7,3) -- (0,3) -- (0,0);
\foreach \i in {0,...,4}{
    \foreach \j in {0,...,3}
    {
        \pgfmathsetmacro{\cubexl}{\i*0.6  + 3.}
        \pgfmathsetmacro{\cubeyl}{\j*0.6  + 0.5}
        \pgfmathsetmacro{\dcubex}{0.2}
        \pgfmathsetmacro{\dcubey}{0.2}
        \draw[draw=black,fill=BurntOrange!40!white] (\cubexl,\cubeyl) -- ++(\dcubex,0) -- ++(0,\dcubey) -- ++(-\dcubex,0) -- cycle;
    }
    \fill [pattern = north east lines] (-0.6,0) rectangle (-0.3, 3.0);
}
\draw (-0.15, 0.15) circle (0.15);
\draw (-0.15, 3 - 0.15) circle (0.15);
\draw (-0.3, 1.5 - 0.17320508075688773) -- (0, 1.5) -- (-0.3, 1.5 + 0.17320508075688773) -- (-0.3, 1.5 - 0.17320508075688773);

\draw[thick,->] (0.0, 1.5) -- (2.5, 1.5) node[anchor=north]{$x$};
\draw[thick,->] (0.0, 1.5) -- (0.0, 4.0) node[anchor=south]{$y$};

\draw[thick,->] (8.5, 1.5) -- (7.5, 1.5) node[anchor=north]{$F_x$};
\draw[thick,->] (8.5, 1.5) -- (8.5, 2.5) node[anchor=south]{$F_y$};

\end{tikzpicture}
    \caption{Schematic of the fiber (orange) reinforced thin plate.}
\label{FIG: PLATE}
\end{figure}
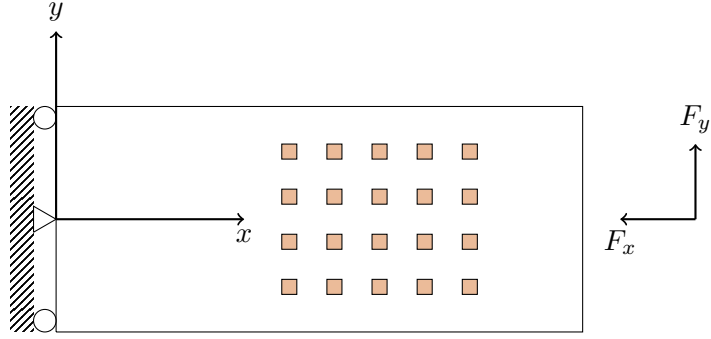

Consider first a thin linear elastic fiber reinforced rectangular plate $[0,L]\times[-c, c]$ with width $L = 100$ and height $2c = 20$. The plate is supported on the left edge $x = 0$, a pinned support at the center point and vertical roller supports at both corners, and subjected to a distributed load along the right edge $x = L$. (See \cref{FIG: PLATE}). The distributed load on the right edge is
\begin{equation}
F =  \Big(-\frac{3  p L} {2  c^2} , \frac{3  p (1 - (y / c)^2)} {4  c} \Big),\, y \in [-c, c]
\label{EQ:PLATE_BC}
\end{equation}
here $p$ is a load strength parameter. Both the matrix and the reinforcing fibre are made of homogeneous and isotropic elastic materials for which the Young's moduli and Poisson ratios are listed in \cref{TAB: PLATE}.
This is a multiscale composite material problem, generally, to resolve each fiber is computationally unaffordable.
Therefore, the homogenized constitutive relation will be applied. \added{The homogenized solutions are used as reference solutions in this section.}

\begin{table}[htbp]
\begin{center}
\begin{tabular}{ l c c }
\toprule
   Materials & Young's modulus  & Poisson ratio  \\
\midrule
Matrix & 1000 & 0.49999 \\
 Fiber & 3000 & 0.39999  \\
\bottomrule
\end{tabular}
\caption{Material properties of fiber reinforce thin plate.}
\label{TAB: PLATE}
\end{center}
\end{table}

\deleted{Using mathematical homogenization,} The governing linear elastostatics equations with plane stress assumptions are expressed in terms of the (Cauchy) stress components $\sigma_{ij}$,
\begin{equation}
\begin{aligned}
    \sigma_{ij, j} + b_i &= 0 \textrm{ in } \Omega\\
    u_i &= \bar{u}_i \textrm{ on } \Gamma_{u}\\
    n_j \sigma_{ji} &= \bar{t}_i \textrm{ on } \Gamma_{t}\\
\end{aligned}
\label{EQ: PLATE_DF}
\end{equation}
here $u_i$ is the displacement, $\Omega$,  $\Gamma_{u}$, and $\Gamma_{t}$ are the computational domain, the displacement boundary, and the traction boundary. Summation convention is employed for repeated indices. The strain tensor is
\begin{equation}
\varepsilon_{mn} = \frac{1}{2}\Big(\frac{\partial u_n}{\partial x_m} +  \frac{\partial u_m}{\partial x_n}\Big)
\end{equation}
The \added{homogenized} linear constitutive relation between strain and stress is written as
\begin{equation}
\sigma_{ij} = \overline{\mathsf{C}}_{ijmn} \varepsilon_{mn}
\end{equation}
here $\overline{\mathsf{C}}_{ijmn} $ are the homogenized constitutive tensor components. Corresponding to our learning problem, we have
\begin{equation*}
  \boxed{\begin{aligned}
	\mathcal{M}_{\bt}(u) &= \overline{\mathsf{C}}_{ijmn} \varepsilon_{mn}(u)\\
	\bt &= \overline{\mathsf{C}}_{ijmn} 
\end{aligned}}
\end{equation*}

The computational domain is discretized by $24 \times 12$ quad elements with linear shape functions. The solution has the finite element approximation $\bu^h = \bv^h + \bar{\bu}^h$, where $\bv^h$ represents unknowns and $\bar{\bu}^h$ represents the boundary states.
With any test function  $\bw^h$, the integration form of \cref{EQ: PLATE_DF} is written as

\begin{equation}
\begin{aligned}
\int_{\Omega}\varepsilon(\bw^h)^T \overline{\mathsf{C}}\varepsilon(\bv^h) d\Omega = \int_{\Omega} \bw^{h}_i  b_i d\Omega + 
\added{ \int_{\Gamma_{t}} \bw^{h}_i  \bar{t}_i d\Gamma_{t} }
- \int_{\Omega}\varepsilon(\bw^h) \overline{\mathsf{C}} \varepsilon(\bar{\bu}^h) d\Omega
\end{aligned}
\label{EQ: PLATE_INT}
\end{equation}

Integrating \cref{EQ: PLATE_INT} with the 2 point Gaussian quadrature in each direction, the fully discretized governing equation becomes
\begin{equation*}
\bK(\overline{\mathsf{C}}) \bv^h - \bF = \bm0
\end{equation*}
here $\bK$ is the stiff matrix depends on the homogenized constitutive tensor $\overline{\mathsf{C}}$, and $\bF$ is the external force vector.

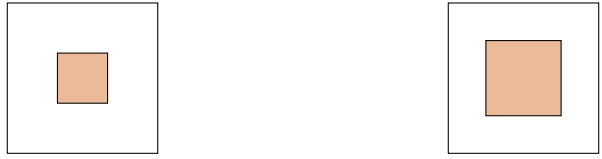
\begin{figure}[htbp]
\centering
\begin{tikzpicture}
  \pgfmathsetmacro{\dcube}{2}
  \draw (0,0) --++(\dcube,0) --++(0,\dcube) -- ++(-\dcube,0) -- (0,0);
  \pgfmathsetmacro{\dcubef}{0.6666666666}
  \draw[draw=black,fill=BurntOrange!40!white] (\dcubef,\dcubef) -- ++(\dcubef,0) -- ++(0,\dcubef) -- ++(-\dcubef,0) -- cycle;
\end{tikzpicture}~~~~~~~~~~~~~~~~~~~~~~~~~~~~~
\begin{tikzpicture}
  \pgfmathsetmacro{\dcube}{2}
  \draw (0,0) --++(\dcube,0) --++(0,\dcube) -- ++(-\dcube,0) -- (0,0);
  \pgfmathsetmacro{\dcubef}{0.5}
  \draw[draw=black,fill=BurntOrange!40!white] (\dcubef,\dcubef) -- ++(2*\dcubef,0) -- ++(0,2*\dcubef) -- ++(-2*\dcubef,0) -- cycle;
\end{tikzpicture}
    \caption{Schematic of the fine-scale unit-cell problems with fiber volume fraction 1/9 (left) and 1/4 (right).}
\label{FIG: PLATE UNIT CELL}
\end{figure}

Several solution data pairs $(\bv_k, p_k), k=1,...,N$ are collected by varying the load strength $p$ in \cref{EQ:PLATE_BC} applied on the right edge.
The constitutive relation used to generate data is obtained through the homogenization procedure discussed in~\cite{feyel2000fe2, yuan2008toward}.
For each pair $mn = 11, 22, \textrm{ or } 12$, a fine-scale unit-cell problem (See \cref{FIG: PLATE UNIT CELL}) resolving micro-scale features on $\Omega^f$ is constructed,
\begin{equation}
\begin{aligned}
    \mathsf{C}_{ijkl}(\varepsilon^{f,mn}_{kl} +  I_{klmn})_{,j} &= 0 \textrm{ in } \Omega^f\\
    u^{f,mn}_i(y) &= u^{f,mn}_i(y + Y) \textrm{ on } \partial \Omega^f\\
    u_i^{f,mn}(y) &= 0 \textrm{ on } \partial \Omega^{f,\,vert}
\end{aligned}
\label{EQ: PLATE_FINE_DF}
\end{equation}
with the superposition of a background strain $I_{klmn} =(\delta_{mk}\delta _{nl} + \delta_{nk}\delta_{ml})/2$, due to the macro-scale strain, and the correction displacement $u^{f,mn}$.  The correction displacement $u^{f,mn}$ is assumed to be periodic in all directions, and zero at all corners $\partial \Omega^{f,\,vert}$ of the unit cell. Here $\mathsf{C}_{ijkl}$ denotes the constitutive tensor,  which depends on the material properties in \cref{TAB: PLATE}, and $\varepsilon^{f,mn}$ denotes the strain associated to the correction displacement $u^{f,mn}$ in the unit-cell problem. By solving the fine-scale unit-cell problem \cref{EQ: PLATE_FINE_DF}, the homogenized constitutive tensor component $\overline{\mathsf{C}}_{ijmn}$ is given as
\begin{equation}
\overline{\mathsf{C}}_{ijmn} = \frac{1}{\Omega^f}\int_{\Omega^f}  \sigma^{f,mn}_{ij} d\Omega^f
\end{equation}
where $\sigma^{f,mn}_{ij} = \mathsf{C}_{ijkl}(\varepsilon^{f,mn}_{kl} +  I_{klmn})$.

For the linear constitutive relation, it is unnecessary to use a neural network for approximation. Instead, we only need to learn the entries of a symmetric matrix. However, the algorithm remains the same and automatic differentiation is the workhorse for the optimization. In all the experiments below, we minimize the loss function \cref{eq:loss} using the \texttt{L-BFGS-B} optimizer. The maximum iteration is 5000 and the tolerance for the gradients norm and the relative change in the objective function is $10^{-12}$.

For the UQ analysis, the only QoI is chosen to be the maximum principal stress
\begin{equation}
\sigma_1 = \frac{\sigma_{11} + \sigma_{22}}{2} + \sqrt{\Big(\frac{\sigma_{11} - \sigma_{22}}{2}\Big)^2 + \sigma_{12}^{2}}
\end{equation}
Its gradient forms the basis of the reduced subspace for $\Delta \overline{\mathsf{C}} $.

\subsubsection{Space-invariant Constitutive Relation}

In this case, the volume fraction of the fiber is assumed to be constant of $1/9$ for the whole plate. The fine-scale unit cell problem \cref{EQ: PLATE_FINE_DF}  is solved, according to \cite{yuan2008toward}, to build the homogenized constitutive tensor component $\overline{\mathsf{C}}_{ijmn}$ as follows
\begin{align}
    \overline{\mathsf{C}} = 
    \begin{bmatrix}
    \overline{C}_{1111}    &     \overline{C}_{1122}    &     \overline{C}_{1112}\\
    \overline{C}_{1122}    &     \overline{C}_{2222}    &     \overline{C}_{2212}\\
    \overline{C}_{1112}    &     \overline{C}_{2212}    &     \overline{C}_{1212}\\
    \end{bmatrix}=
    \begin{bmatrix}
    1491.24     &   701.024  &     0\\
    701.024     &   1450.24  &     0\\
    0           &   0        &     362.941
    \end{bmatrix}
\end{align}
One data point ($N = 1$) is generated with load strength $p = 20$, with the linear homogenized multiscale constitutive relation. For the inverse problem, the loss function \cref{eq:loss} is minimized to reach an optimal value of $10^{-12}$ within 50 steps. The following constitutive tensor is obtained
\begin{align}
\overline{\mathsf{C}}_{\bt} = & \begin{bmatrix}
    1491.24     &    701.024    &   -1.12339\times 10^{-8}\\
    701.024     &    1450.24    &   -1.3416\times 10^{-8}\\
   -1.12339\times 10^{-8}  &  -1.3416\times 10^{-8}  & 362.941
\end{bmatrix}
\end{align}
The constitutive relation is recovered exactly from the proposed learning process for the linear case.

\subsubsection{Space-varying Constitutive Relation}\label{sect:nonuniform}
\begin{figure}[htbp]
\centering
\includegraphics[width=0.8\textwidth,keepaspectratio]{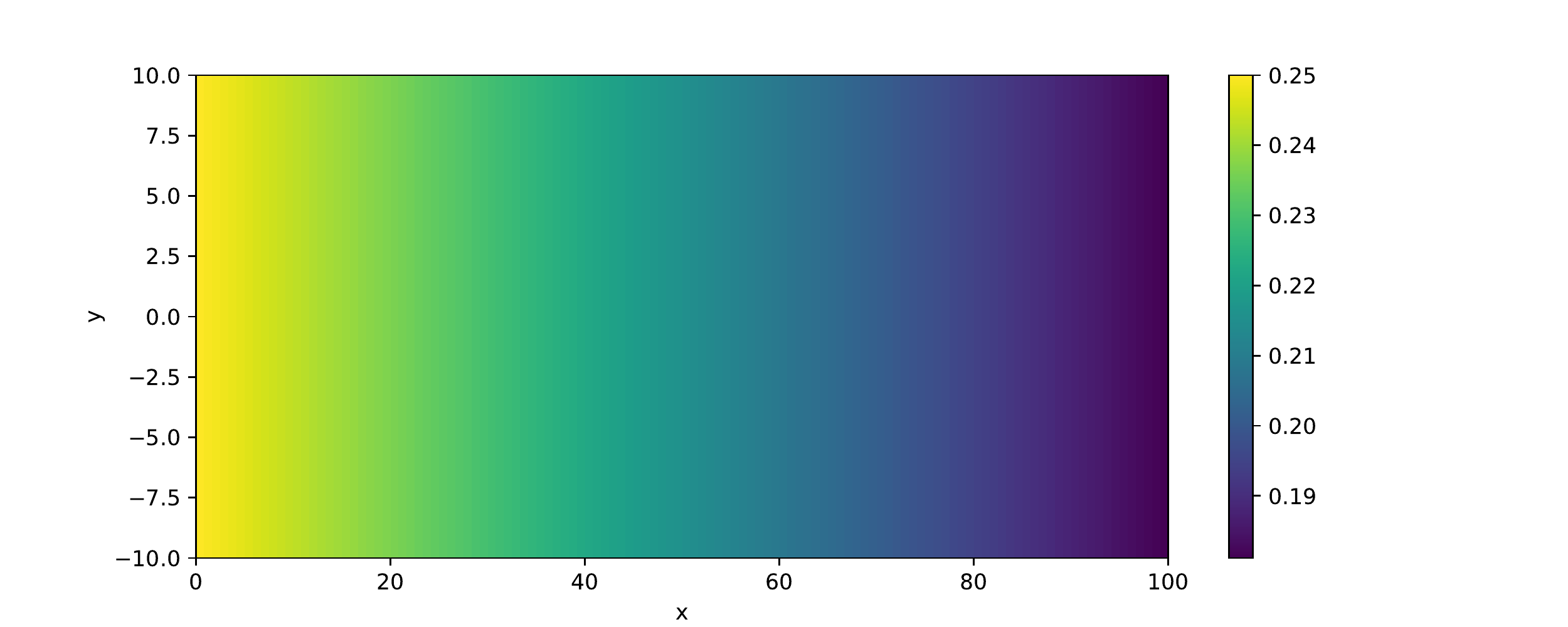}

\includegraphics[width=0.8\textwidth,keepaspectratio]{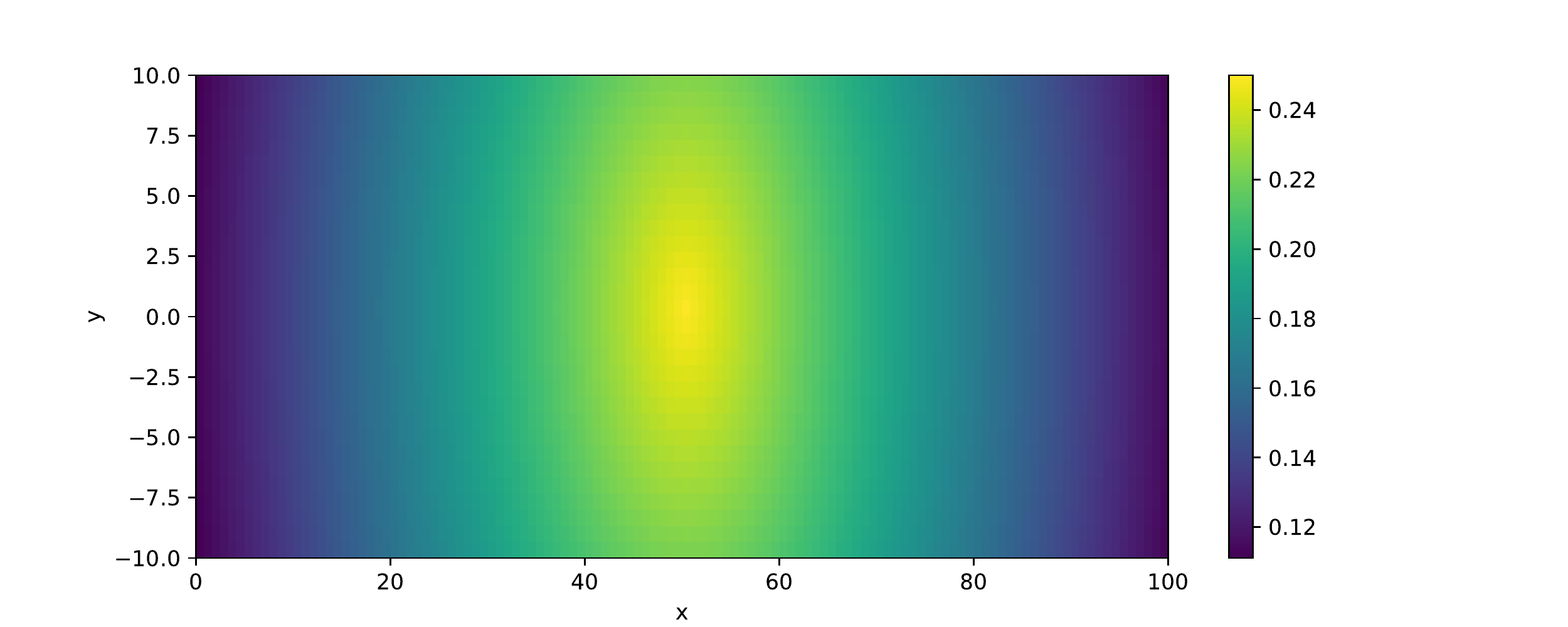}
\caption{The volume fraction distributions of fibers in the thin plane with space-varying constitutive relations in \Cref{sect:nonuniform}.}
\label{fig:PlaneStress_VolumeFraction_12}
\end{figure}

For this case, the volume fraction of the fiber is assumed to be space-varying, between  $1/9$ and $1/4$. Consider two fiber volume fraction distributions depicted in \cref{fig:PlaneStress_VolumeFraction_12}, the fibers volume fraction decreases from the left to the right (\cref{fig:PlaneStress_VolumeFraction_12}-top) as follows,
\begin{equation}
\frac{1}{9}\frac{x}{2L} + \frac{1}{4}\left(1 - \frac{x}{2L}\right)
\label{eq:volume_fraction_1}
\end{equation}
and the fibers are denser at the center of the plate and gradually become sparse along the radius direction (\cref{fig:PlaneStress_VolumeFraction_12}-bottom), as follows
\begin{equation}
\frac{1}{9}\sqrt{\frac{(x - L/2)^2 + y^2}{(L/2)^2 + c^2}} + \frac{1}{4}\Big(1 - \sqrt{\frac{(x - L/2)^2 + y^2}{(L/2)^2 + c^2}}\Big)
\label{eq:volume_fraction_2}
\end{equation}

The homogenized constitutive tensor component $\overline{\mathsf{C}}'_{ijmn}$ obtained by the fine-scale unit
cell corresponding to volume fraction $1/4$ (see \Cref{FIG: PLATE UNIT CELL}-right) is
\[
    \overline{\mathsf{C}}' = \begin{bmatrix}
        1695.92    &      747.42    &     0\\
  747.42     &   1633.96     &     0\\
    0   &  0&  405.76
    \end{bmatrix}
\]
The homogenized constitutive tensor at each element for volume fraction between $\frac{1}{9}$ and $\frac{1}{4}$ is linearly interpolated between $\overline{\mathsf{C}}$ and $\overline{\mathsf{C}}'$.
Therefore, at each element, the constitutive relations are different. Our data-driven algorithm can homogenize the material
based on the global response, and compute the homogenized constitutive tensor.
The learned constitutive relations for these two fiber volume fraction distributions from the observation with load strength $p=20$ are
\begin{equation}
    \mathsf{C}_{\bt1} = \begin{bmatrix}
        1582.58  &    698.793    &   1.24528\\
  698.793  &   1512.1     &    2.80921\\
    1.24528  &    2.80921 &  377.979
    \end{bmatrix} \quad \textrm{ and } \quad
    \mathsf{C}_{\bt2} = \begin{bmatrix}
         1673.94  &    738.872     & 2.59714\\
  738.872 &   1578.87     &  6.06215\\
    2.59714  &   6.06215  & 399.123
    \end{bmatrix}
\end{equation}

These learned constitutive relations $\mathsf{C}_{\bt1}$ and $\mathsf{C}_{\bt2}$ are averages between $\overline{\mathsf{C}}$ and $\overline{\mathsf{C}}'$ that minimize \cref{eq:loss}. For the UQ analysis, the variance of the reduced coordinate of the first case is $\hat \bSig_{\blam} \approx  1.15\times 10^5$,  while in the second case, the variance is $\hat \bSig_{\blam} \approx  4.21\times 10^5$.

\begin{figure}[htbp] 
\centering
\includegraphics[width=0.45\textwidth,keepaspectratio]{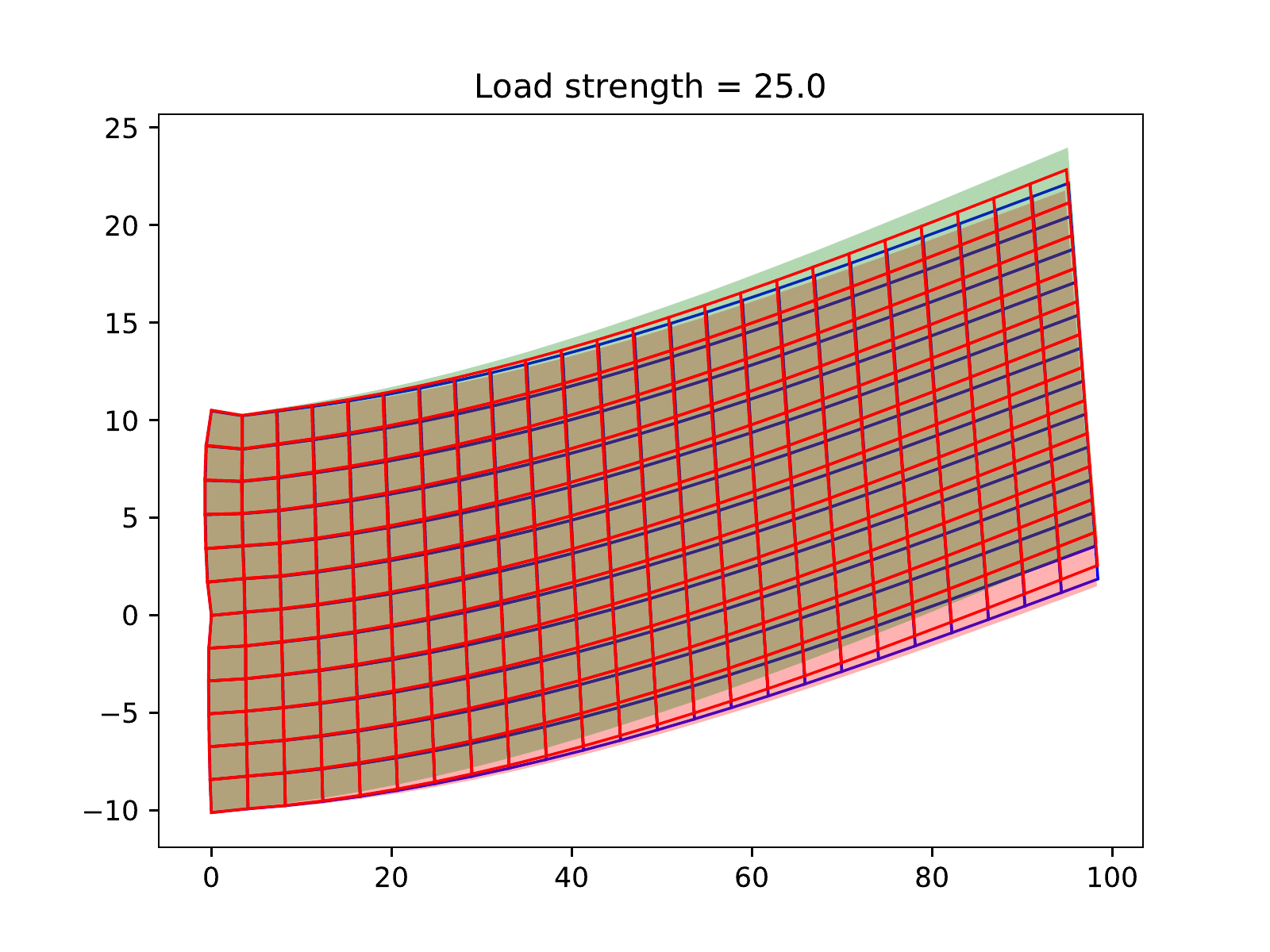}~
\includegraphics[width=0.45\textwidth,keepaspectratio]{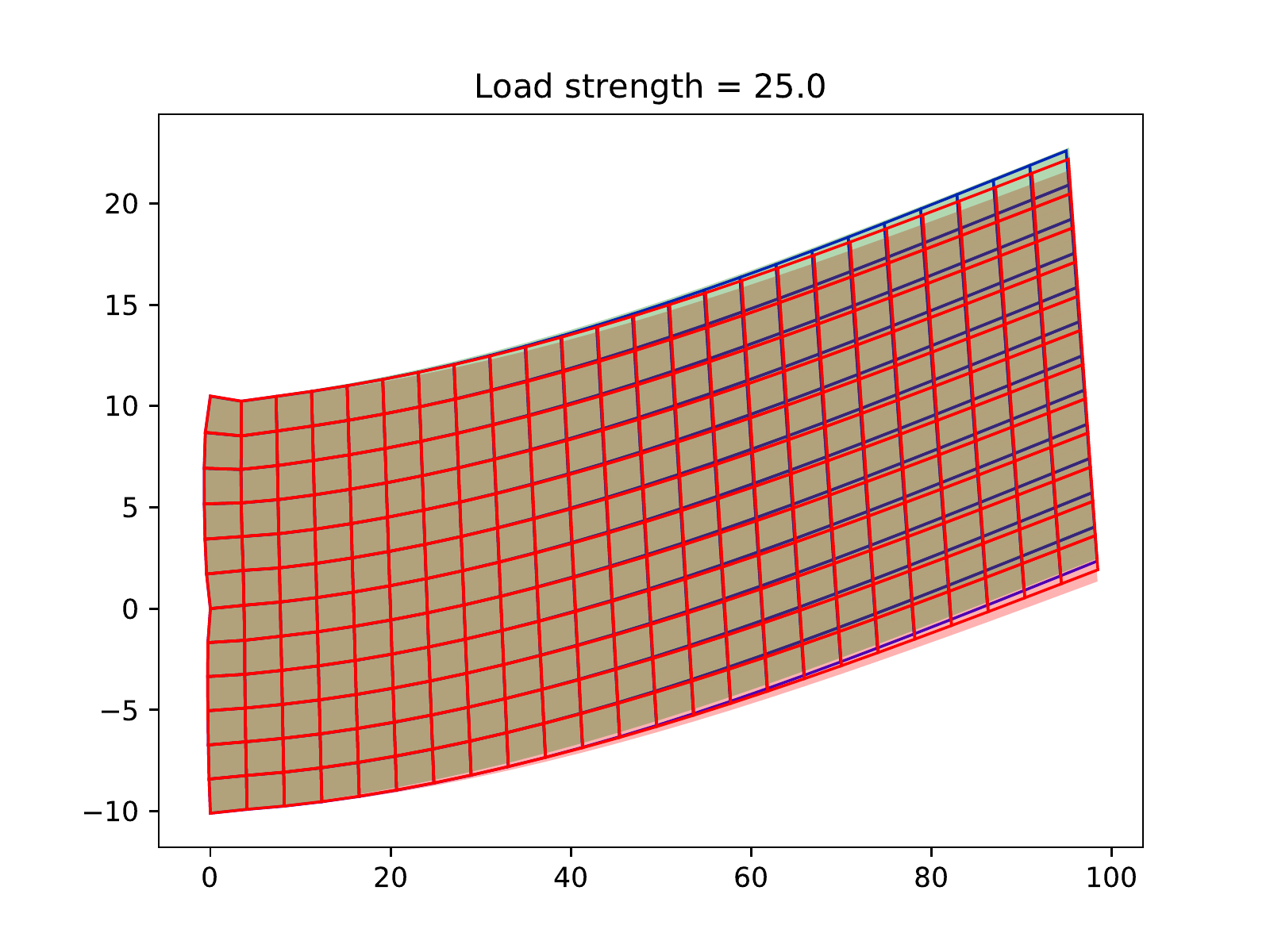}\\
\caption{\added{Predicted shape of thin plate with space-varying fiber distributions} \Cref{eq:volume_fraction_1}~(left) and  \Cref{eq:volume_fraction_2}~(right). \deleted{The grey mesh is the undeformed configuration,} The blue mesh is the exact deformed configuration. The red mesh is the predicted deformed configuration. The green and the pink regions correspond to the quantiles 0.95 and 0.05 results.}
\label{fig:PlaneStress_noise1}
\end{figure}

\begin{figure}[htbp] 
\centering
\scalebox{0.8}{
\begin{tikzpicture}

\definecolor{color1}{rgb}{0.12156862745098,0.466666666666667,0.705882352941177}
\definecolor{color0}{rgb}{1,0.647058823529412,0}

\begin{axis}[
legend cell align={left},
legend style={at={(0.03,0.97)}, anchor=north west, draw=white!80.0!black},
tick align=outside,
tick pos=left,
x grid style={white!69.01960784313725!black},
xlabel={Load strength},
xmin=20.4, xmax=77.6,
xtick style={color=black},
y grid style={white!69.01960784313725!black},
ylabel={$y$ Displacement},
ymin=0.261603788809439, ymax=33.5089132809266,
ytick style={color=black}
]
\path [draw=green!50.19607843137255!black, ultra thick]
(axis cs:25,1.77284512936022)
--(axis cs:25,3.94259677202325);

\path [draw=green!50.19607843137255!black, ultra thick]
(axis cs:35,6.43151725109562)
--(axis cs:35,9.49698690006289);

\path [draw=green!50.19607843137255!black, ultra thick]
(axis cs:45,11.1318952652807)
--(axis cs:45,15.0556859265932);

\path [draw=green!50.19607843137255!black, ultra thick]
(axis cs:55,15.7547424183103)
--(axis cs:55,20.6858921546303);

\path [draw=green!50.19607843137255!black, ultra thick]
(axis cs:65,20.5379792833246)
--(axis cs:65,26.1428034751984);

\path [draw=green!50.19607843137255!black, ultra thick]
(axis cs:75,25.4574147736601)
--(axis cs:75,31.9976719403759);

\addplot [semithick, red, mark=*, mark size=3, mark options={solid}, only marks]
table {%
25 2.84344005787245
35 7.98081608102158
45 13.1181921041705
55 18.2555681273195
65 23.3929441504683
75 28.5303201736175
};
\addlegendentry{Prediction}
\addplot [semithick, blue, mark=triangle*, mark size=3, mark options={solid}, only marks]
table {%
25 2.14832435471132
35 7.00765409659584
45 11.8669838384803
55 16.7263135803649
65 21.5856433222494
75 26.4449730641338
};
\addlegendentry{Reference}
\end{axis}

\end{tikzpicture}}~~~~~
\scalebox{0.8}{
\begin{tikzpicture}

\definecolor{color0}{rgb}{0.12156862745098,0.466666666666667,0.705882352941177}
\definecolor{color1}{rgb}{1,0.647058823529412,0}

\begin{axis}[
legend cell align={left},
legend style={at={(0.03,0.97)}, anchor=north west, draw=white!80.0!black},
tick align=outside,
tick pos=left,
x grid style={white!69.01960784313725!black},
xlabel={Load strength},
xmin=20.4, xmax=77.6,
xtick style={color=black},
y grid style={white!69.01960784313725!black},
ylabel={\(\displaystyle y\) Displacement},
ymin=0.207942439235659, ymax=29.6483628746543,
ytick style={color=black}
]
\path [draw=green!50.19607843137255!black, ultra thick]
(axis cs:25,1.54614336811832)
--(axis cs:25,2.77883688227201);

\path [draw=green!50.19607843137255!black, ultra thick]
(axis cs:35,6.24031397814553)
--(axis cs:35,7.87485141829097);

\path [draw=green!50.19607843137255!black, ultra thick]
(axis cs:45,10.8847387410053)
--(axis cs:45,12.966521479154);

\path [draw=green!50.19607843137255!black, ultra thick]
(axis cs:55,15.5219189738136)
--(axis cs:55,18.1054666253203);

\path [draw=green!50.19607843137255!black, ultra thick]
(axis cs:65,20.1196574608889)
--(axis cs:65,23.2666542534746);

\path [draw=green!50.19607843137255!black, ultra thick]
(axis cs:75,24.7822702067159)
--(axis cs:75,28.3101619457716);

\addplot [semithick, red, mark=*, mark size=3, mark options={solid}, only marks]
table {%
25 2.18075223247506
35 7.05305312546509
45 11.9253540184551
55 16.797654911445
65 21.6699558044354
75 26.5422566974248
};
\addlegendentry{Prediction}
\addplot [semithick, blue, mark=triangle*, mark size=3, mark options={solid}, only marks]
table {%
25 2.59769407609375
35 7.63677170653116
45 12.6758493369687
55 17.7149269674063
65 22.7540045978436
75 27.7930822282811
};
\addlegendentry{Reference}
\end{axis}

\end{tikzpicture}}
\caption{\added{Prediction of the $y$ displacement of the upper right corner of the thin plate with space-varying fiber distributions} \Cref{eq:volume_fraction_1}~(left) and \Cref{eq:volume_fraction_2}~(right), and the confidence intervals corresponding to the quantiles 0.95 and 0.05, subjected to different external loads.}
\label{fig:PlaneStress_noise1_bar}
\end{figure}
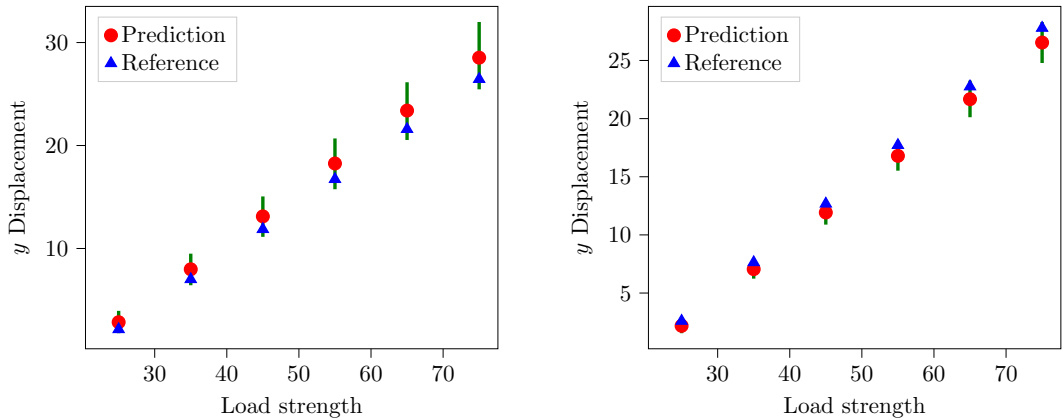

We then apply the learned constitutive relations and its corresponding variance $\hat \bSig_{\blam}$ to predict displacements for various load strengths $(p = 25,\,35,\,45,\,55,\,65,\,75)$ and perform UQ analysis. 
The exact solution, the predicted solution,  and the confidence intervals corresponding to the quantiles 0.95 and 0.05, constructed with 2000 samples generated by \cref{EQ:MC} are depicted in \cref{fig:PlaneStress_noise1,fig:PlaneStress_noise1_bar}. 
The predicted result and the exact solution are in good agreement, the confidence regions of the constitutive relation fluctuations contain exact solutions.

\subsection{Rubber Membrane}
\begin{figure}[htbp]
\centering
   \centering
\begin{tikzpicture}
 \draw[thick,->] (0.0, 0.0) -- (6.0, 0.0) node[anchor=north]{$r$};
 \draw[thick,->] (0.0, 0.0) -- (0.0, 4.0) node[anchor=south]{$z$};
 \draw[BurntOrange, line width=1.0mm] (5,0) arc (30:90:5.773502691896258);
 \draw[BurntOrange, line width=1.0mm, dash pattern=on 6pt off 3pt] (0,0) -- (4.0, 0);
 \draw[thick,->] (4.0, -0.2) -- (5.0, -0.2) node[anchor=north west]{$\bar{u}_r$};

  \draw[thick,->] (1.19543396, 1.96142029) -- (1.34486321, 2.51909782);
  \draw[thick,->] (2.30940108 ,1.5) -- (2.59807621, 2. ) ;
    \draw[thick,->] (3.26598632 ,0.76598632) -- (3.67423461, 1.17423461) ;
 \draw[thick,->] (1.76753906, 1.76721677) -- (1.98848144, 2.30061887);
 \draw[thick,->] (2.8117486,  1.16434212) -- (3.16321717, 1.62238488);

\draw node at (2.2,0.9) {$p$};

\end{tikzpicture}
\caption{Schematic of the axisymmetric rubber membrane: the initial undeformed configuration (dash line) and the current deformed configuration (solid line).}
\label{FIG: MEMBRANE}
\end{figure}
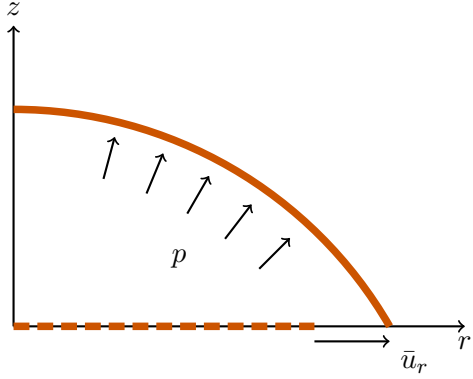

Finally, the rubber membrane with both material and geometric nonlinearities is considered. The circular membrane of radius $L= 1$ is initially flat, and then stretched and inflated under the non-conservative (follower) pressure load (See \cref{FIG: MEMBRANE}). Due to the axisymmetricity, the model problem is reduced to a 1D problem,  described in the $(r, z)$ plane. A point in the undeformed configuration has coordinate $(R,Z) \in [0,1] \times \{0\}$ with line element $dS$, and in the deformed configuration has coordinate $(r, z)$ with line element $ds$.  The displacement vector is defined as
\[
u_r = r -  R,\, u_z = z - Z
\]
Three principal stretch ratios of the membrane at the point are given by
\begin{equation*}
\lambda_1 = \frac{ds}{dS},\, \lambda_2 = \frac{2\pi r}{2\pi R},\, \textrm{ and } \lambda_3 = \frac{t}{T}
\end{equation*}
here $ds = \sqrt{dr^2 + dz^2}$, $dS = \sqrt{dR^2 + dZ^2}$, and $t$ and $T$ are the thickness of the membrane in the deformed configuration and the undeformed configuration. The rubber membrane is assumed to be incompressible,
\[
dSdR T = dsdr t \, \Rightarrow \, \lambda_1 \lambda_2 \lambda_3 = 1
\]
The corresponding first Piola-Kirchhoff stresses are  $P_1$, $P_2$ and $P_3 = 0$.  The weak form of the governing equations is written as
\begin{equation}
\begin{aligned}
    2\pi T \int_{0}^{1} (P_1 \delta \lambda_1 + P_2 \delta \lambda_2) R dS - 2\pi p\int_{0}^{r}\bm{n} \delta \bm{u} rds = 0 \quad \forall \delta \bm{u}
    \end{aligned}
    \label{EQ: MEM_INT}
\end{equation}
Here the first term is the virtual internal work and the second term is the virtual external work corresponding to the non-conservative (follower) pressure load $p$. Due to the symmetry and pre-stretch of the membrane, the boundary conditions are
\begin{equation}
u_r(0,0) = 0,\, u_r(1,0) = \bar{u}_r ,\, u_z(1,0) = 0
\end{equation}

The computational domain is discretized by 100 elements with linear shape functions. Integrating \cref{EQ: MEM_INT} with 3 points Gaussian quadrature in each element lead to the fully discretized governing equation
\begin{equation*}
\bP(\bu , \Mcal) = \bF(\bu, p)
\end{equation*}
here $\bP$ and $\bF$ correspond to the discretization of the first and second term in \cref{EQ: MEM_INT}, and $\Mcal$ represents the parameters $P_1$, $P_2$, which are unknowns and will be calibrated with observations.  Here $\Mcal$ is the constitutive relation related to the principal stretches to the first Piola-Kirchhoff stresses.

Several solution data $(\bu_k, p_k)$ are collected by varying the pressure load $p$. For generating the data, the rubber membrane is a Mooney-Rivlin hyperelastic incompressible material \cite{fried1982finite} with a energy density function $W$ as follows,
\begin{equation}
\begin{aligned}
&W(\lambda_1, \lambda_2, \lambda_3) =  \mu(\lambda_1^2 + \lambda_2^2 + \lambda_3^2 - 3) + \alpha(\lambda_1^2\lambda_2^2 + \lambda_2^2\lambda_3^2 + \lambda_3^2\lambda_1^2 - 3)\\
&J = \lambda_1 \lambda_2 \lambda_3 = 1
\end{aligned}
\label{EQ: Mooney-Rivlin-hyperelastic}
\end{equation}
here $\mu$ and $\alpha$ are model parameters. The true constitutive relation is given as
\begin{equation}
\begin{aligned}
P_1 = \frac{\partial W}{\partial \lambda_1}\qquad \textrm{and} \qquad
P_2 = \frac{\partial W}{\partial \lambda_2}
\end{aligned}
\label{EQ:MEM_CONSTITUTIVE_LAW}
\end{equation}
After nondimensionalizing the problem parameters, $p' = \frac{p}{\mu T}$, $\alpha' = \frac{\alpha}{\mu}$, $P_1' = \frac{P_1}{\mu}$, and $P_2' = \frac{P_2}{\mu}$, the integration form \cref{EQ: MEM_INT} becomes
\begin{equation}
\begin{aligned}
    \int_{0}^{1} (P_1' \delta \lambda_1 + P_2'  \delta \lambda_2) R dS - p'\int_{0}^{r}\bm{n} \delta \bm{u} rds = 0 \quad \forall \delta \bm{u}
    \end{aligned}
    \label{EQ: MEM_INT_NONDIM}
\end{equation}
The constitutive relation is approximated by a neural network
\begin{equation*}
\boxed{\Mcal_{\bt} : (\lambda_1, \lambda_2) \rightarrow (P_1', P_2')}
\end{equation*}
where $\bt$ is the weights and biases of our neural network.

In all the experiments below, we minimize the loss function \cref{eq:loss} using the \texttt{L-BFGS-B} optimizer. The maximum iteration is 20,000 and the tolerance for the gradients norm and the relative change in the objective function is $10^{-12}$.

For the UQ analysis, the only QoI is chosen to be the maximum principal stress, in this case, the maximum principal stress is $\max\{P_1', P_2'\}$. Its gradient forms the basis of the reduced subspace for $\Delta \bt$.

\begin{remark}
In the hyper-elasticity case, $(P_1,P_2)$ are gradients of $W$. In this case, to enforce this relationship, we can model $W$ directly with a neural network and use automatic differentiation to compute $(P_1,P_2)$. In this way, the model is guaranteed to preserve energy because the stress is the gradient of an energy function $W$. \added{We have not pursued this approach in the study since this approach requires second-order derivatives of the neural network when computing the gradients of \mbox{\Cref{EQ:MEM_CONSTITUTIVE_LAW}}.}
\end{remark}

\subsubsection{Space-invariant Constitutive Relation}\label{sect:ucl}

Consider the Mooney-Rivlin hyperelastic incompressible rubber membrane with parameter $\alpha' = 0.1$, subjected to non-dimensional pressure loads $p'$ of 0.0, 0.5, 1.0, \ldots, 7.5, 8.0. We collect 17 data points $(\bu_k, p_k)$ to train the FEM-neural network framework. The neural network consists of two hidden layers with 20 neurons in each layer\replaced{. We use \texttt{tanh} as activation functions since it is smooth compared to other choices, such as ReLU or leaky ReLU. \texttt{tanh} is also bounded, which prevents creating extremely large or small values in the intermediate layers.}{; the activation function is \texttt{tanh} since it is smooth}. \added{In the numerical examples below, } the neural network converges within 15,000 iterations.

\begin{figure}[htbp] 
\centering
\includegraphics[width=0.5\textwidth,keepaspectratio]{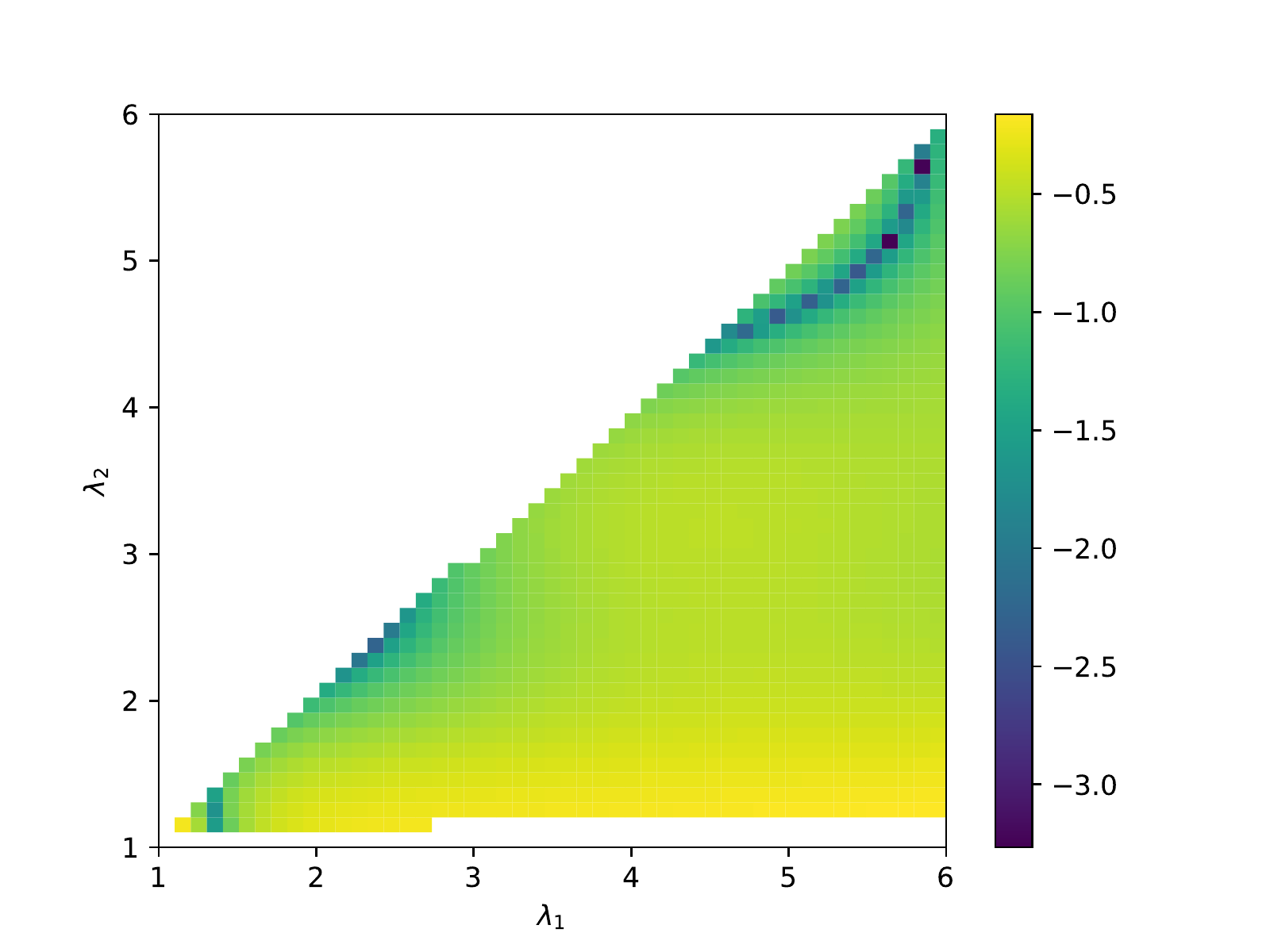}~
\includegraphics[width=0.5\textwidth,keepaspectratio]{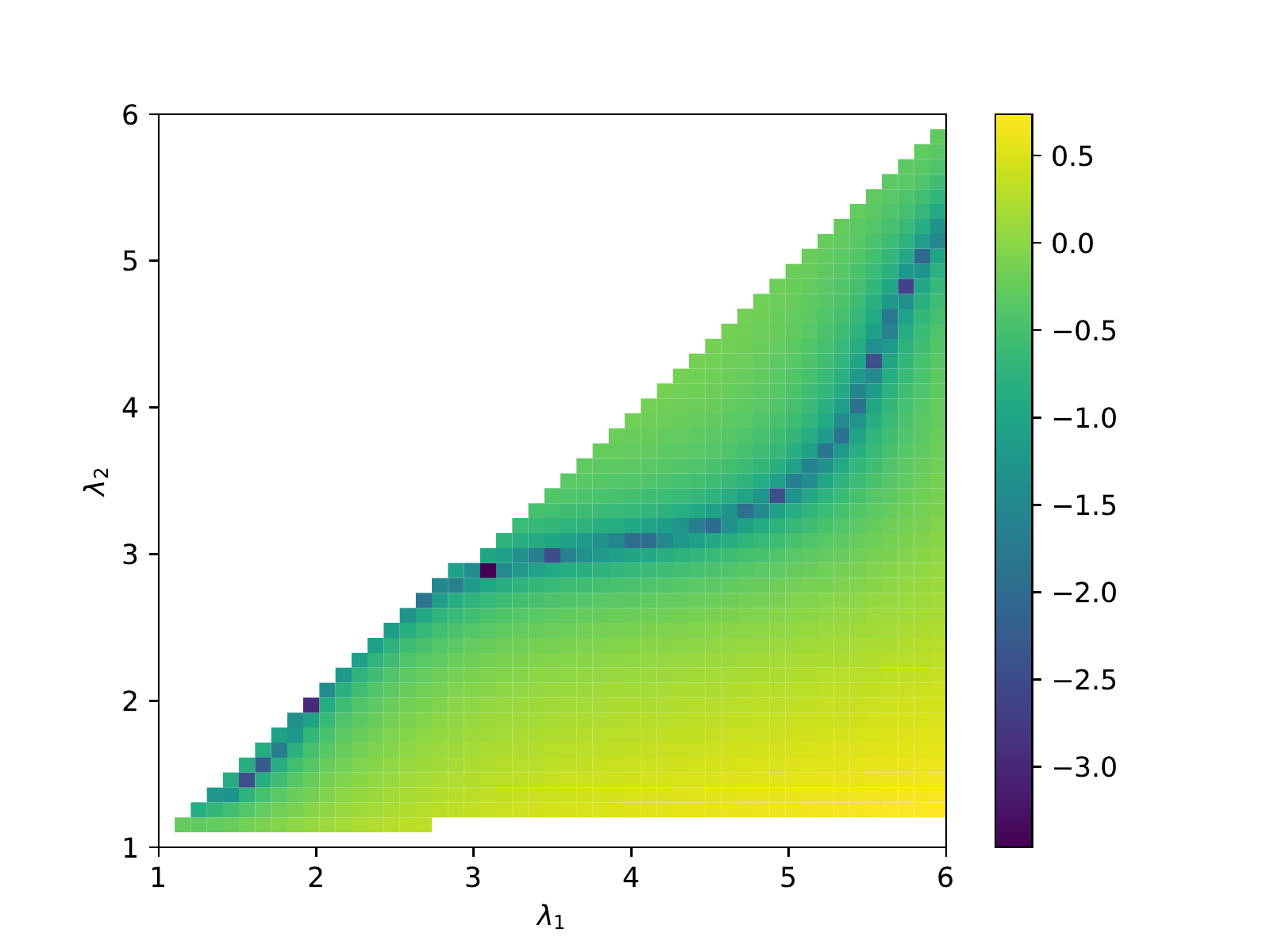}

\caption{The logarithmic (base 10) absolute error of the surrogate constitutive relation learned by the neural network approach for the rubber membrane in the space-invariant constitutive relation case, the first component of the stress $P_1'$ (left), the second component of the stress $P_2'$ (right). The choice of $(\lambda_1,\lambda_2)$ pairs used for this figure corresponds to the convex hull of the $(\lambda_1,\lambda_2)$ pairs found in the training data. This is why not all points are present in these plots.}
\label{fig:lambda_v1}
\end{figure}

\Cref{fig:lambda_v1} shows the differences between the calibrated and the exact constitutive relations. Since the present supervised learning framework is incapable of extrapolating, predicting the corresponding stresses for stretches $\lambda_1$, $\lambda_2$ that are outside of the convex hull of the given 17 data points, we do not expect the neural network result to be accurate for all $(\lambda_1, \lambda_2)$. Therefore, in the figure, only a subset of strains $\lambda_1, \lambda_2$ that appears in the inputs is depicted. For the first component of the stress, $P_1'$ with respect to the strain $(\lambda_1, \lambda_2)$, the calibrated constitutive relation almost overlaps with the exact one. For the second component $P_2'$, the calibrated result deviates a little from the exact one near $\lambda_1=6$, $\lambda_2=1$. This is because the loss function, the residual force \cref{eq:loss}, is less sensitive to the second component compared with the first component. However, we have treated both components equally in the formulation of the loss function. The backpropagation will update the first component more {effectively} than the second one.

We then apply the learned constitutive relation and its corresponding variance $\hat \bSig_{\blam} \approx 2.43\times 10^{-6}$ to predict the displacements for various load strengths ($p' = 2.2,\,4.2,\,6.2$) and perform UQ analysis. \Cref{fig:v31} shows the exact displacements, the predicted displacements and the confidence intervals corresponding to the quantiles 0.95 and 0.05, constructed with 2000 samples generated by \cref{EQ:MC} in the $z$ and $r$ directions. \deleted[id=r2]{We need to point out we are also able to extrapolate to predict the behavior of the rubber membrane subjected to a pressure load of $8.2$, although the neural network has not seen any observations beyond pressure $8.0$.}

The exact displacements and the predicted displacements almost collapse on top of each other. As for the quantification of the homogenized error,  the blue uncertainty region is almost vanished, which means the predicted value is of high reliability. This is consistent with the fact that there is no noise in the underlying constitutive model.

\begin{figure}[htbp]
\centering
\scalebox{0.8}{\input{figures/v31_2.2_1__}}~
\scalebox{0.8}{\input{figures/v31_2.2_2__}}
\scalebox{0.8}{\input{figures/v31_4.2_1__}}~
\scalebox{0.8}{\input{figures/v31_4.2_2__}}
\scalebox{0.8}{\input{figures/v31_6.2_1__}}~
\scalebox{0.8}{\input{figures/v31_6.2_2__}}
\caption{The space-invariant constitutive relation case: the predicted (dash red line) and the exact displacements (blue line) in the $z$ and $r$ directions (left to right) of the rubber membrane, subjected to pressure loads $2.2$, $4.2$, and $6.2$\deleted[id=r2]{, and $8.2$ (top to bottom),} and their corresponding confidence intervals (blue region) corresponding to the quantiles 0.95 and 0.05.}
\label{fig:v31}
\end{figure}


\subsubsection{Space-varying Constitutive Relation}\label{sect:noisycons}

\begin{figure}[htbp]
\centering
\scalebox{1.0}{
\begin{tikzpicture}

\definecolor{color0}{rgb}{0.12156862745098,0.466666666666667,0.705882352941177}

\begin{axis}[
tick align=outside,
tick pos=left,
x grid style={white!69.01960784313725!black},
xlabel={R},
xmin=-0.05, xmax=1.05,
xtick={-0.2,0,0.2,0.4,0.6,0.8,1,1.2},
xticklabels={,0.0,0.2,0.4,0.6,0.8,1.0,},
y grid style={white!69.01960784313725!black},
ytick={0.090,  0.095, 0.100, 0.105, 0.110},
yticklabels={0.090,  0.095, 0.100, 0.105, 0.110},
ylabel={$\alpha^{'}$},
ymin=0.09, ymax=0.11
]
\addplot [semithick, color0, forget plot]
table [row sep=\\]{%
0	0.098 \\
0.333333333333333	0.108 \\
0.666666666666667	0.093 \\
1	0.106 \\
};




\end{axis}

\end{tikzpicture}}
\caption{The parameter $\alpha'$ in the energy density function \cref{EQ: Mooney-Rivlin-hyperelastic} at different elements used  in the space-varying constitutive relation case. The curve is generated by connecting randomly sampled points $\alpha'(0)$, $\alpha'(\frac{1}{3})$, $\alpha'(\frac{2}{3})$, and  $\alpha'(1)$  with $10\%$ error around 0.1.}
\label{fig:Membrane_alpha}
\end{figure}
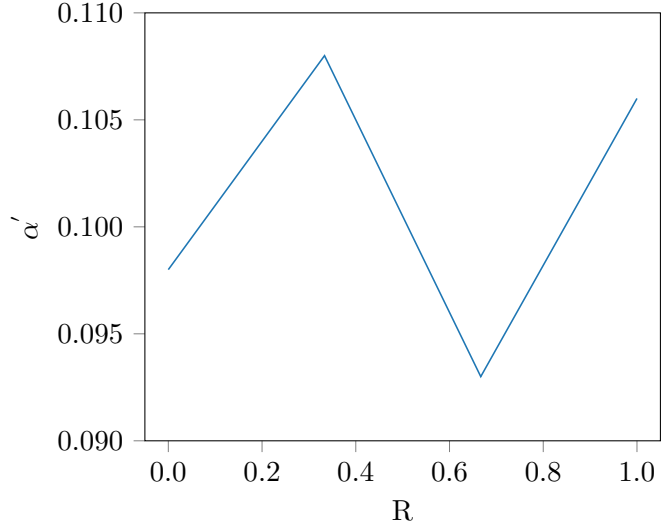

Due to the thickness and material variations of the rubber membrane, the constitutive relation can be space-varying or noisy.  The noise is incorporated in the test problem by varying the parameter $\alpha'$ in the energy density function \cref{EQ: Mooney-Rivlin-hyperelastic} at different elements. The $\alpha'$ used to generate data is randomly sampled, which represents the noise in the constitutive relation due to the thickness and material variations of the rubber membrane. The $\alpha'(R)$ curve is depicted in  \cref{fig:Membrane_alpha}, which is assumed to be continuous, and about 10\% variation is included.

We generate 17 data with the same non-dimensional pressures as in the previous space-invariant constitutive relation case. The displacements of the underlying model are solved with heterogeneous constitutive relations.
Besides, we adopt the same neural network structure, optimizer and training iterations as that in the space-invariant case.

\begin{figure}[htbp] 
\centering
\includegraphics[width=0.5\textwidth,keepaspectratio]{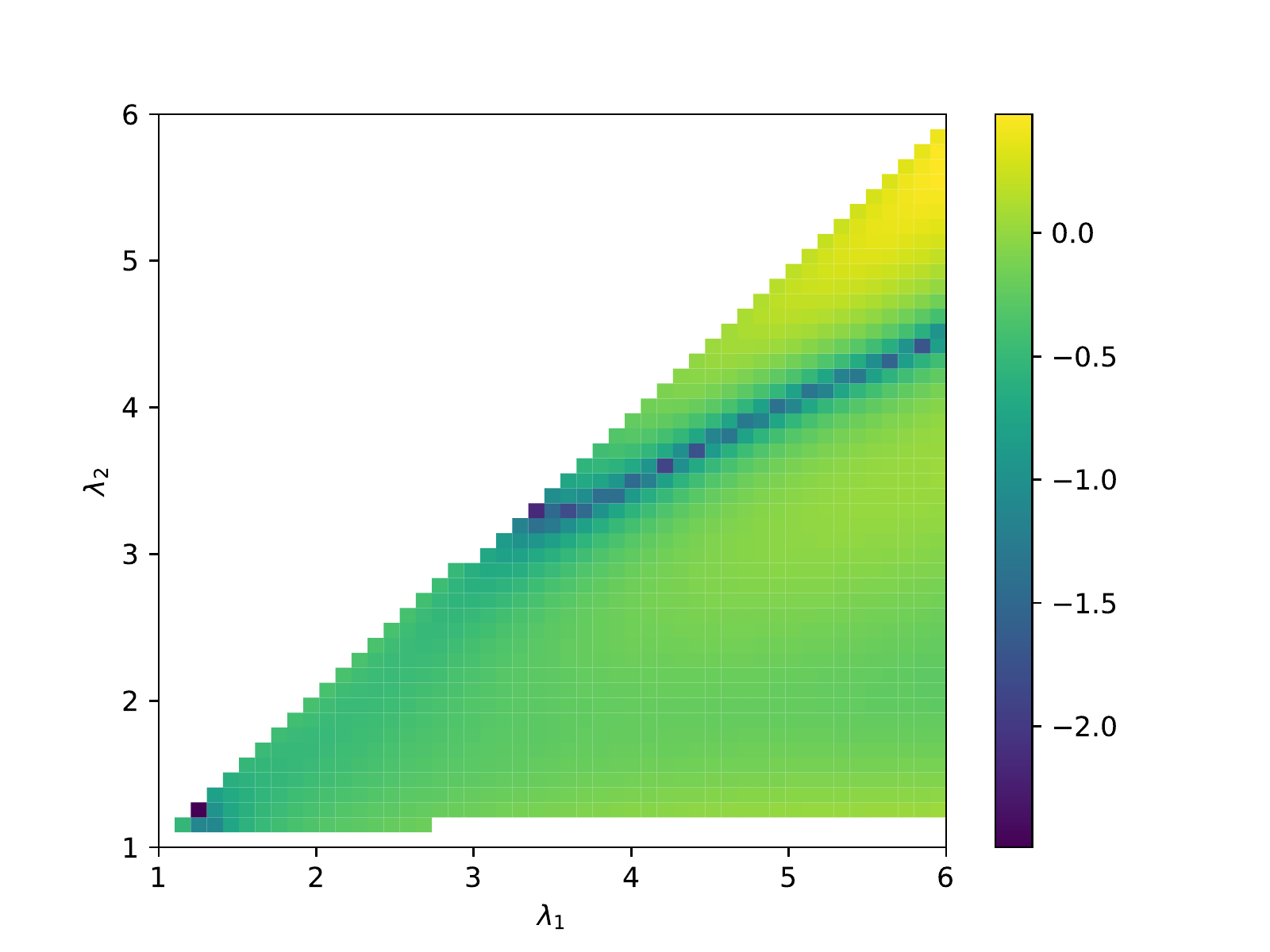}~
\includegraphics[width=0.5\textwidth,keepaspectratio]{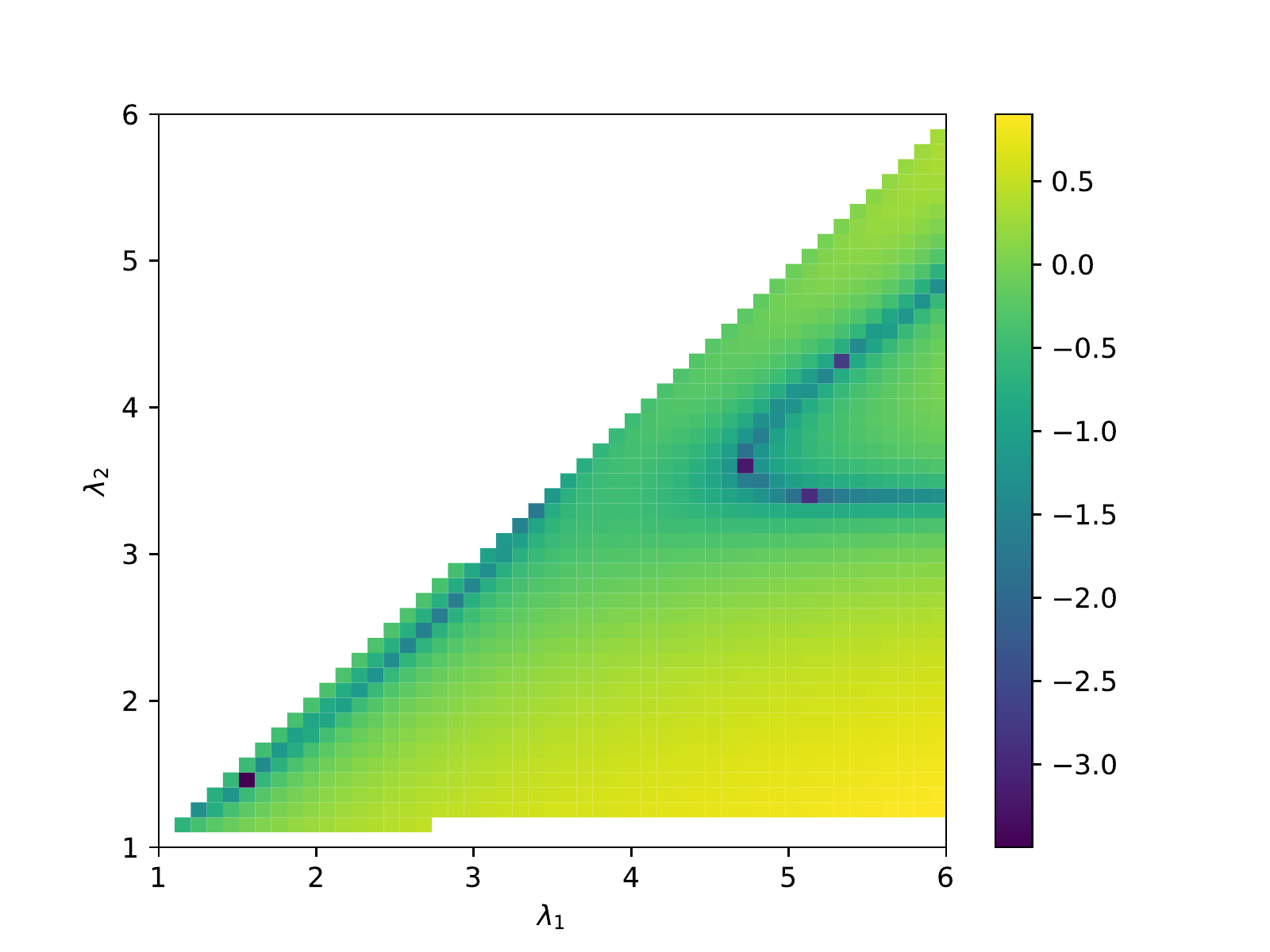}
\caption{The logarithmic (base 10) absolute error of the surrogate constitutive relation learned by the neural network approach for the rubber membrane in the space-varying constitutive relation case, the first component of the stress $P_1'$~(left), the second component of the stress $P_2'$~(right). }
\label{fig:lambda_v1_3}
\end{figure}

\Cref{fig:lambda_v1_3} shows the difference between the calibrated relation and the reference constitutive relation corresponding to $\alpha' = 0.1$. Note in this case, the homogenized constitutive relation~(reference) is not necessarily the true relation but only serves as a reference. What is more relevant is how we can predict the behavior of the rubber membrane subjected to the external pressure load. Since the training data cover the range $p\in[0, 8.0]$, the test pressures are chosen to be $p' = 2.2,\, 4.2,\, \textrm{ and } 6.2$, which avoids extrapolations. The predicted displacements and exact displacements, obtained with the heterogeneous constitutive relation are shown in \cref{fig:v31_2}. Good agreement can be observed. The estimated variance of the reduced coordinate is $\hat \bSig_{\blam} \approx  2.23\times 10^{-6}$. The confidence intervals corresponding to the quantiles 0.95 and 0.05 are constructed with 2000 samples generated by \cref{EQ:MC}. The exact solution lies in the uncertainty region, which demonstrates the effectiveness of the presented framework.

\begin{figure}[htbp]
\centering
\includegraphics[width=0.5\textwidth,keepaspectratio]{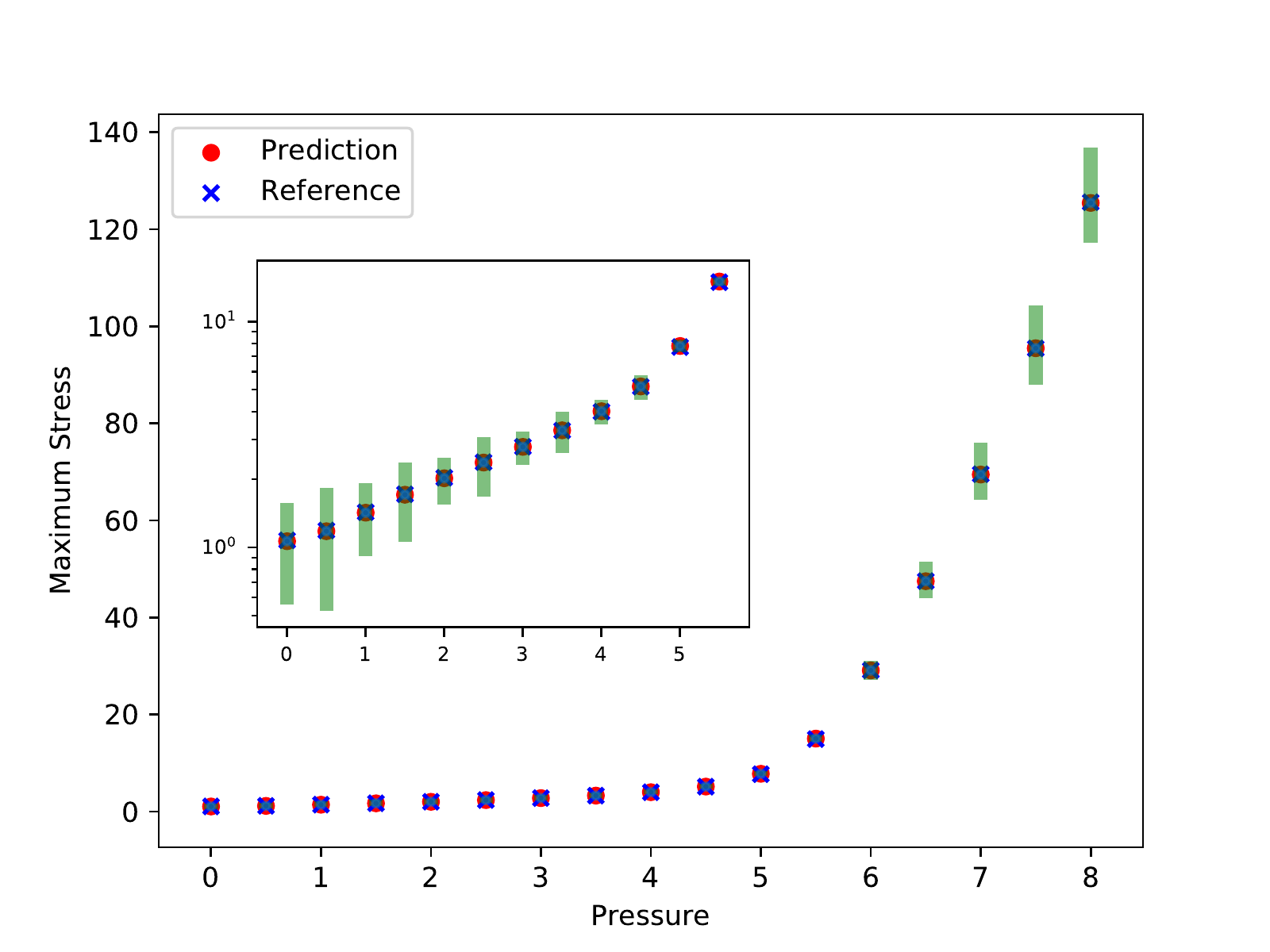}~
\includegraphics[width=0.5\textwidth,keepaspectratio]{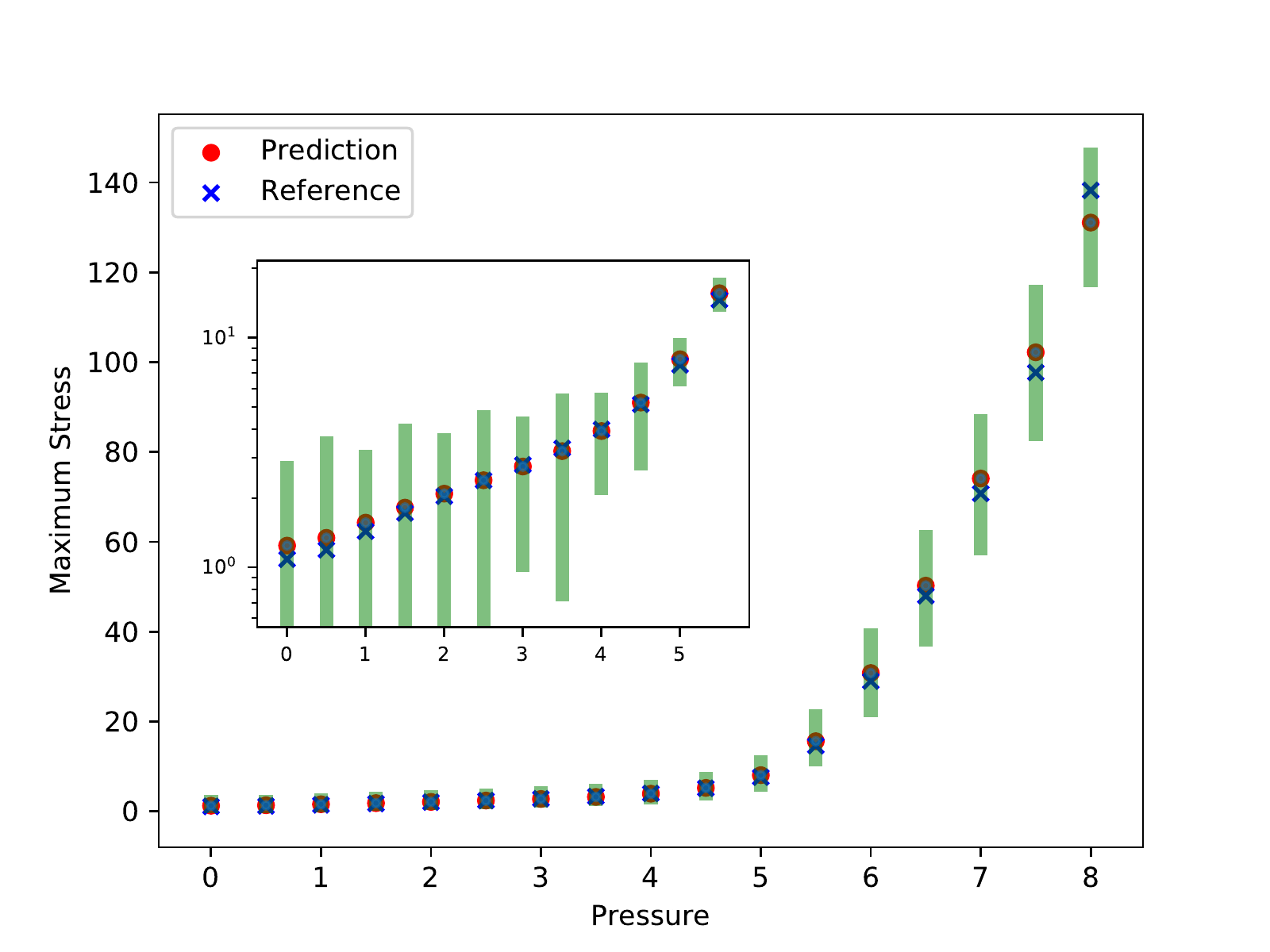}
\caption{\added{Prediction of the maximum stresses of the rubber membrane and their $99.7\%$ confidence intervals, subjected to various pressure loads:  the space-invariant constitutive relation case (left) and the space-varying constitutive relation case (right). Plots with logarithmic scale for y-axis are shown separately for clarity.}}
\label{fig:stress}
\end{figure}

It is also interesting to visualize the uncertainties for the maximum stress, we have used as the QoI. \Cref{fig:stress} shows the maximum stresses together with its $3\delta$ confidence interval for both the space-invariant constitutive relation case and space-varying constitutive relation case. We can see that the exact maximum stresses lie in the confidence interval in all cases. For the space-varying constitutive relation case, we have larger uncertainty error bounds, which is consistent with our intuition since the space-invariant constitutive relation case does not have homogenization errors.


\begin{figure}[htbp]
\centering
\scalebox{0.8}{\input{figures/v33_2.2_1__}}~
\scalebox{0.8}{\input{figures/v33_2.2_2__}}
\scalebox{0.8}{\input{figures/v33_4.2_1__}}~
\scalebox{0.8}{\input{figures/v33_4.2_2__}}
\scalebox{0.8}{\input{figures/v33_6.2_1__}}~
\scalebox{0.8}{\input{figures/v33_6.2_2__}}
\caption{The space-varying constitutive relation case: the predicted (dashed red line) and the exact displacements (blue line) in the $z$ and $r$ directions (left to right) of the rubber membrane, subjected to pressure loads $2.2$, $4.2$, and $6.2$ (top to bottom), and their corresponding confidence intervals (blue region) corresponding to the quantiles 0.95 and 0.05.}
\label{fig:v31_2}
\end{figure}

\subsection{Comparison with Other Approximations}
 \label{sect:cplrs_rbf}
\begin{figure}[htbp]
\centering
\includegraphics[width=0.5\textwidth,keepaspectratio]{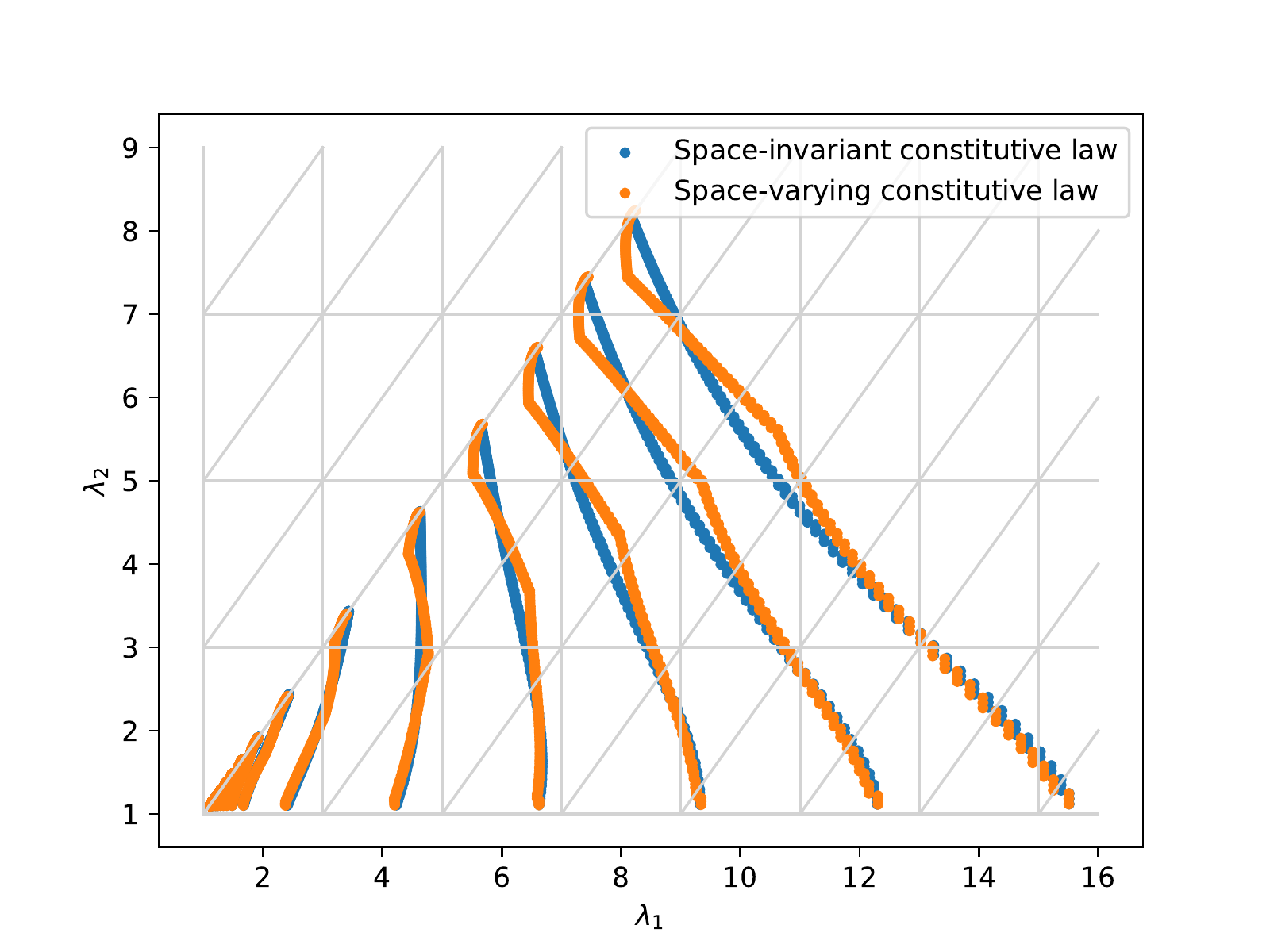}
\caption{The principal stretch pairs $(\lambda_1$, $\lambda_2)$ that appear in the training set for the space-invariant constitutive relation case and the space-varying constitutive relation case and the triangulation of the parametric domain with $h = 2.0$.}
\label{fig:membrane_strain_paris}
\end{figure}

The constitutive relation \cref{EQ:MEM_CONSTITUTIVE_LAW} can also be approximated by other functions, besides neural networks. 
In this section, the proposed neural network approach is compared with piecewise linear functions~(PL) and radial basis functions~(RBF). 
These functions approximate the constitutive relations  on the parameter space, which is $[0,20]^2$ and the superscript 2 corresponds to the dimension of the input variables
The parameter space is chosen based on \cref{fig:membrane_strain_paris}, which depicts all the principal stretch pairs $(\lambda_1$, $\lambda_2)$ post-processed from the training set for both space-invariant constitutive relation and space-varying constitutive relation cases.  
In practice, it is impossible to adopt a point-to-point surface fitting since the stress data is not directly available from the final output. 
Although in our case the strain distribution from the training data can be obtained, it can be arbitrarily nonuniform and lie on some low-dimensional manifold~(See \cref{fig:membrane_strain_paris}), 
which challenges both PL and RBF approximations.
However, we show that neural networks lend us a universal method for approximating the constitutive relation without any prior information on the data distribution. Besides, it shows robustness in terms of noise and generalizes better than fine-tuned PL and RBF. Consequently, the neural network approach has the potential to achieve state-of-the-art results using data-driven modeling.

For comparison, we use \texttt{L-BFGS-B} optimizer for all cases. The optimization is terminated when the objective function is called 15000 times. The training data consists of 17 samples as mentioned before. The test data consists of three samples subjected to non-dimensional pressure loads of $p=2.2$, $4.2$, $6.2$. The neural network (NN) is the standard fully-connected deep nets with two 20-neuron hidden layers. And losses evaluated on both the training set and the test set at each training iteration are reported.

\subsubsection{Comparison with Piecewise Linear Functions}\label{sect:cplrs}

For the piecewise linear functions~(PL-$h$), we use a uniform triangulation on $[0,20]^2$ with mesh size $h =$ 0.4, 1.0, 2.0, which  are chosen based on~\cref{fig:membrane_strain_paris}. As for the partitioning, a uniform grid is chosen~(See \cref{fig:membrane_strain_paris} for $h=2.0$). Adaptive triangulation (e.g., Delaunay triangulation) is possible but difficult to implement robustly.  The number of parameters is shown in \cref{tab:param}.

\begin{table}[htbp]
\centering
\begin{tabular}{@{}lllll@{}}
\toprule
Model & PL-0.4  & PL-1.0 & PL-2.0 & NN  \\
\midrule
Degrees of freedom          & 5000 & 400 & 100 & 520\\
\bottomrule
\end{tabular}
\caption{Number of parameters for different surrogate functions.}
\label{tab:param}
\end{table}

\begin{figure}[htbp]
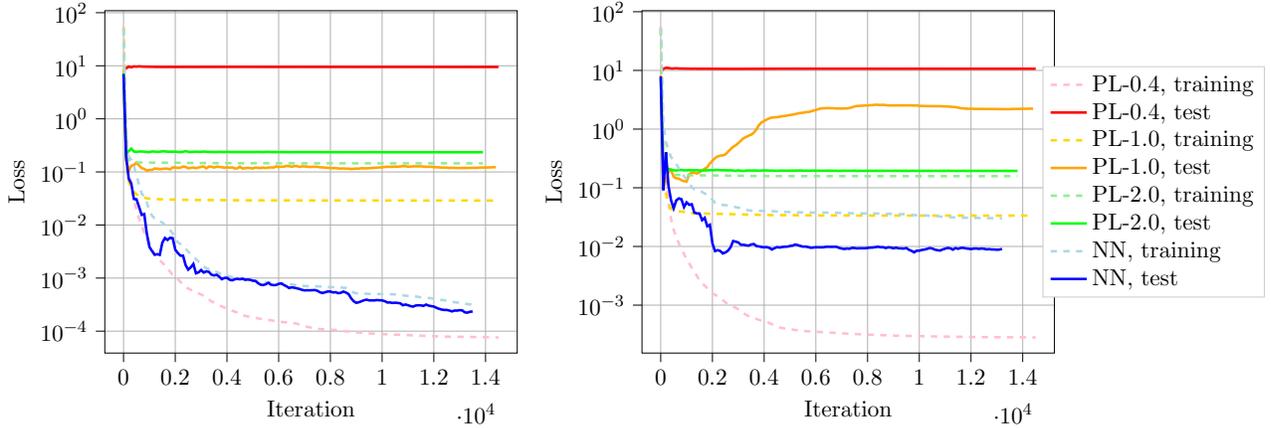

\centering
\scalebox{0.8}{\input{figures/loss1}}~
\scalebox{0.8}{\input{figures/loss3}}
\caption{The losses evaluated on both the training set and the test set at each training iteration with PL-0.4, PL-1.0 and PL-2.0 and NN for the space-invariant constitutive relation case (left) and the space-varying constitutive relation case (right).}
\label{fig:losses}
\end{figure}

\begin{figure}[htbp] 
\centering
\begin{subfigure}[b]{0.5\linewidth}
\centering
\includegraphics[width = 1.0\textwidth, keepaspectratio]{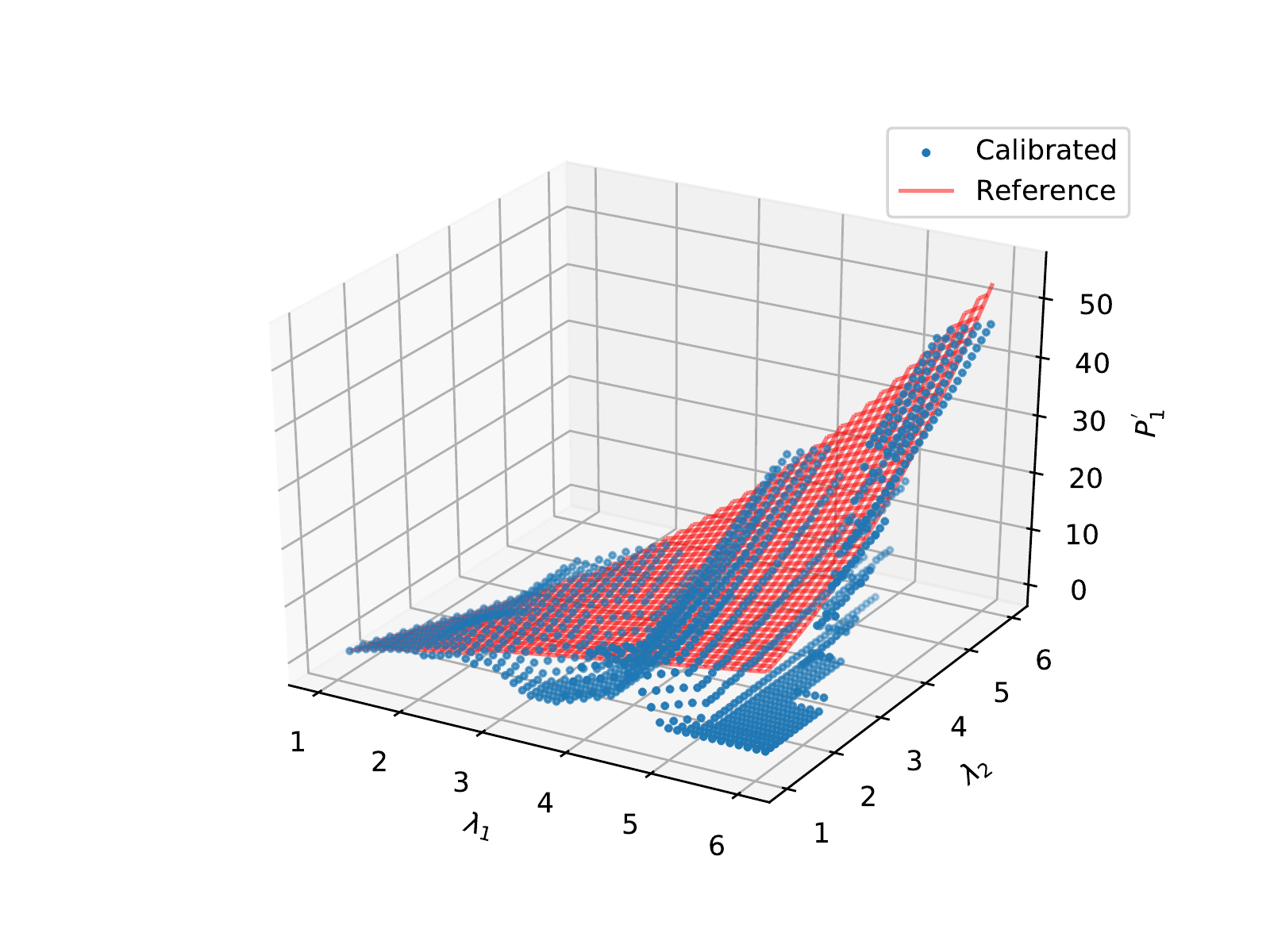}
\caption{PL-0.4}
\end{subfigure}~
\begin{subfigure}[b]{0.5\linewidth}
\centering
\includegraphics[width = 1.0\textwidth, keepaspectratio]{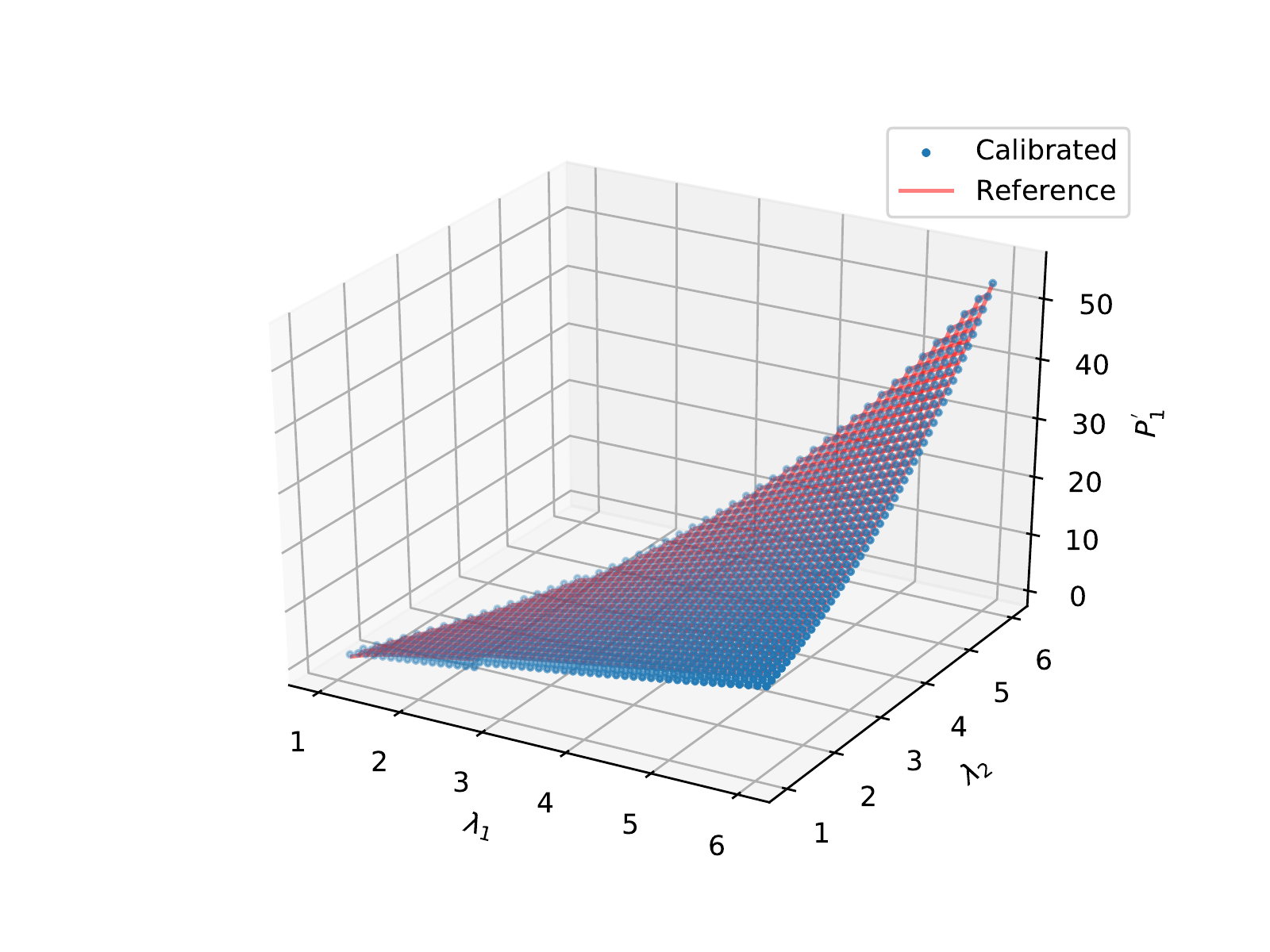}
\caption{NN}
\end{subfigure}
\caption{\deleted{Sampled calibrated constitutive relations~(first component of the stress $P_1'$) learned by PL-0.4, PL-1.0, PL-2.0 and NN for the rubber membrane in the space-invariant constitutive relation case. Zero initial guess for PL-$h$ is used. (e) shows the value of $P_1'$ at fixed $\lambda_2=3.0$ for PL-1.0, PL-2.0 and NN. PL-0.4 clearly fails to fit the reference surface.}
\added{Sampled calibrated constitutive relations~(first component of the stress $P_1'$) learned by PL-0.4~(left) and NN~(right) for the rubber membrane in the space-invariant constitutive relation case. Zero initial guess for PL-$h$ is used.}}
\label{fig:l1}
\end{figure}

\Cref{fig:losses} shows losses~\cref{eq:loss} at each training iteration, including the losses evaluated on both the training set and the test set. PL-0.4 achieves the least training loss, thanks to the large number of local degrees of freedom, which will fit the unknown functions on these training data in \cref{fig:membrane_strain_paris}. However, PL-0.4 fails to predict or approximate the stresses associated with principal stretch pairs $(\lambda_1$, $\lambda_2)$ that do not appear in the training data. These stresses fail to be updated during the training process. Therefore, the approximated surfaces are quite rough and highly oscillating~(See \cref{fig:l1}-left). Certainly, the overfitting leads to poor performance on the test set. PL-1.0 and PL-2.0 reach plateaus rapidly, which indicates under-fitting, i.e. failing to capture the underlying trend of the data. Hence, they also lead to poor performance in the test set. And it is worth mentioning  PL-1.0 reaches similar training error in the space-varying case as NN, but struggles to generalize well on the test data.

NN achieves consistent losses on both the training set and the test set. And the losses are reduced by about four orders of magnitude during the training process. We believe NN can regularize the surrogate function and consistently interpolate on the missing input states, hence leads to smooth constitutive relations~(See \cref{fig:l1}-right). And NN can be more efficient for high dimensional problems, compared with polynomial approximations, for which the number of parameters depends on the input parameter dimension exponentially.

\begin{figure}[htbp]
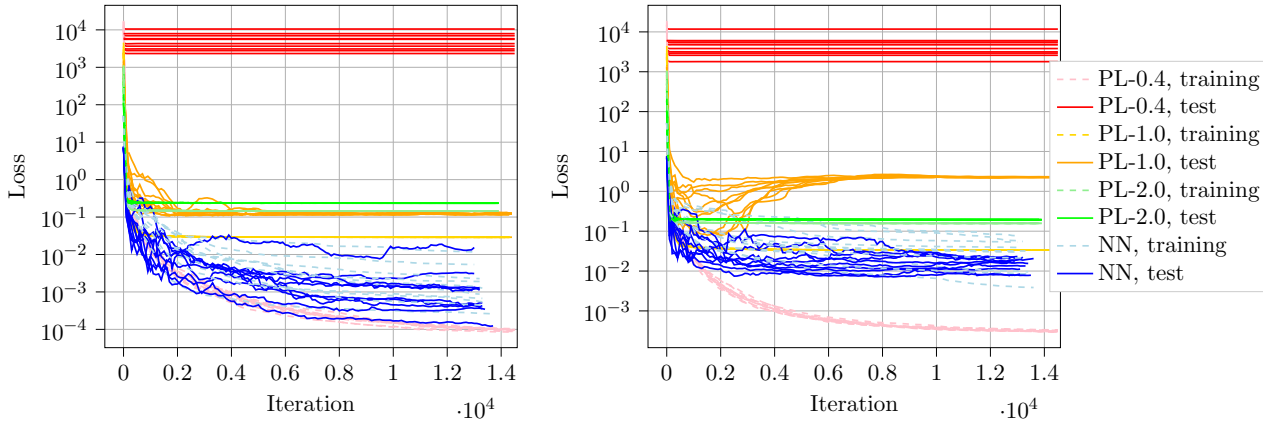

\centering
\scalebox{0.8}{\input{figures/loss1u-initial}}~
\scalebox{0.8}{\input{figures/loss3n-initial}}
\caption{The losses evaluated on both the training set and the test set at each training step with PL-0.4, PL-1.0, PL-2.0 and NN for the space-invariant constitutive relation case~(left) and the space-varying constitutive relation case~(right). Different curves correspond to different initial guesses.}
\label{fig:losses-initial-conditions}
\end{figure}

Moreover,  we also study the sensitivity to initial weights~\cite{2015arXiv151106422M} of these two approaches, since both approaches involve highly non-convex optimization problems. We start from several initial guesses, the i.i.d. Gaussian random variables with mean 0 and standard deviation 100, for both PL and NN approaches~(We also tried standard deviation 10 and found no substantial difference.). \Cref{fig:losses-initial-conditions} depicts that PL is quite robust with respect to different initial guesses for all loss curves overlap well on each other. However, all these initial guesses lead to poor performance on the test set. Meanwhile, the NN approach shows good generalization property, i.e. consistent losses on both training set and test set, even in the space-varying constitutive relation case.

\subsubsection{Comparison with Radial Basis Functions}\label{sect:crbf}

The function approximation combines radial basis functions with global low order polynomials, as follows,
\begin{equation}
    f_{\bt}(\bx) = \sum_{i=1}^n \alpha_i g_\sigma(\|\bx-\bx_i \|) + a + \mathbf{b}^T\bx
\end{equation}
where the parameter $\bt=\begin{pmatrix} \alpha_1, \alpha_2, \ldots, \alpha_n, a, \mathbf{b} \end{pmatrix}$ is to be calibrated through minimizing \cref{EQ: OBJECTIVE}. A good choice of the basis function $g_\sigma$ is important for the approximation. In this paper, we consider the commonly used inverse multi-quadric radial basis function
\begin{equation}
    g_\sigma(r) = \frac{1}{\sqrt{r^2+\sigma^2}}
\end{equation}
where $\sigma$ is the scalar parameter and is chosen to be $\sigma = 20$, since the computational domain is $[0,20]^2$. The centers of the radial basis functions are distributed uniformly in the domain. We consider four cases~(RBF-$h$) with $h=0.1$, $0.4$, $1.0$, $2.0$, which stands for the distance between centers in both directions. The number of parameters is shown in \cref{tab:param2}.

\begin{table}[htbp]
\centering
\begin{tabular}{@{}llllll@{}}
\toprule
Model & RBF-0.1  &RBF-0.4 &RBF-1.0 & RBF-2.0 & NN  \\
\midrule
Degrees of freedom  &  4000       & 2500 & 400 & 100 & 520\\
\bottomrule
\end{tabular}
\caption{Number of parameters for different surrogate functions.}
\label{tab:param2}
\end{table}

\Cref{fig:rbflosses} shows the losses~\cref{eq:loss} at each training iteration,
including the losses evaluated on both the training set and the test set, where NN performs consistently better than RBF-$h$. In the space-invariant case, as we increase the degrees of freedom for RBF-$h$, we achieve better test accuracy. However, this does not hold for the space-varying case.
\begin{figure}[htbp]
\centering
\scalebox{0.8}{\input{figures/rbfloss1.tex}}~
\scalebox{0.8}{\input{figures/rbfloss3.tex}}
\caption{The losses evaluated on both the training set and the test set at each training step with \replaced`{radial basis function}{response surface approaches}~(RBF-0.1, RBF-0.4, RBF-1.0 and RBF-2.0)  and the neural network approach~(NN) for the space-invariant constitutive relation case (left) and the space-varying constitutive relation case (right).}
\label{fig:rbflosses}
\end{figure}

\subsubsection{\added[id=r2]{Comparison with Radial Basis Function Networks}}
\label{sect:rbfn}

\added[id=r2]{Finally, we consider the comparison with radial basis function networks (RBFN) \cite{orr1996introduction}}
$$
f_{\bt}(\bx) = \sum_{i=1}^n w_i g_{\sigma_i}(\|\bx-\bx_i\|)
$$
\added[id=r2]{In this case, the weights, shape parameters, and centers of radial basis functions are optimized with the gradient descent method}
$$\bt=\left(w_1, w_2, \ldots, w_n, \sigma_1, \sigma_2, \ldots, \sigma_n, \bx_1, \bx_2, \ldots, \bx_n\right)$$
\added[id=r2]{The radial basis function is the commonly used Gaussian function}
$$g_{\sigma}(\|\bx\|) = \exp(-\sigma\|\bx\|^2)$$
\added[id=r2]{In the computation, the initial guess for $w_i$ are all zeros, $\sigma_i$ are all ones and $\bx_i$ are uniformly drawn from $[0,7]^2$ (a reasonable computational domain). For numerical stability, we constrain $\sigma_i$ to be nonnegative during the training. We choose different $n=$ 100, 400, 1600, 2500 (the corresponding model is denoted as RBFN-$n$); we compare them with the same neural network architecture in the last section and the number of parameters are shown in} \Cref{tab:param3}.

\begin{table}[htbp]
\centering
\begin{tabular}{@{}llllll@{}}
\toprule
Model & RBFN-100  &RBFN-400 &RBFN-1600 & RBFN-2500 & NN  \\
\midrule
Degrees of freedom  &  400       & 1600 & 6400 & 10000 & 520\\
\bottomrule
\end{tabular}
\caption{Number of parameters for different surrogate functions.}
\label{tab:param3}
\end{table}

\begin{figure}[htbp]
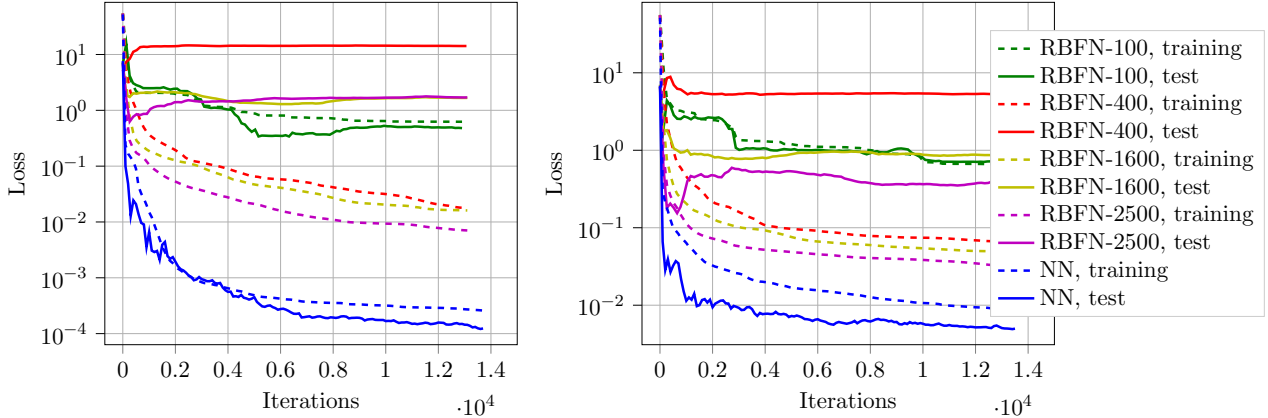

\centering
\scalebox{0.8}{\input{figures/rbfnloss1.tex}}~
\scalebox{0.8}{\input{figures/rbfnloss3.tex}}
\caption{\added{The losses evaluated on both the training set and the test set at each training step with radial basis function network approaches~(RBFN-100, RBFN-400, RBFN-1500 and RBFN-2500)  and the neural network approach~(NN) for the space-invariant constitutive relation case (left) and the space-varying constitutive relation case (right).}}
\label{fig:rbfnlosses}
\end{figure}

\added[id=r2]{\Cref{fig:rbfnlosses} shows the losses at each training and testing iteration. We can see that the neural network approach performs consistently better than RBFN. Particularly, we see that RBFN easily overfits for large $n$; this is because the constraints on $w_i$, $\sigma_i$, $\bx_i$ are weak and RBFN easily fits to the noise. For small $n$ RBFN does not overfit but has poor accuracy.}

\section{Conclusion}\label{sect:conc}

Data-driven approaches continue to gain popularity for constructing coarse-grained models when data from high-fidelity simulations and high-resolution experiments become richer. In this work, a general framework combining traditional FEM and neural networks for predictive modeling is presented. The proposed framework discretizes the physical system through FEM, but replaces the constitutive relation, or any coarse-grained model part, with a black-box neural network. The neural network is trained based on the global response information, for example, the displacement field, instead of point-to-point strain vs stress data; such data are easier to measure from experiments.

Furthermore, the applicability and accuracy of the present framework are analyzed. When the data set contains comprehensive constitutive data, and the relation is smooth enough, the error in the constitutive relation predicted by the neural network is bounded by the optimization error and the discretization error. The uncertainties due to the heterogeneity of the material are quantified efficiently using our FEM-neural network framework. 

The framework was tested on a multi-scale fiber-reinforced thin plate problem and a highly nonlinear rubbery membrane problem. In the rubbery membrane problem, a comparison between deep neural networks and other function approximations is presented, and the strengths of the neural network approach, i.e., its good regularization and generalization properties, were highlighted.

An interesting area for future work would be to extend the current framework to partially observed data. The examples in the present work use the whole displacement field; however, using only part of the displacement field, for example specifying only the displacement field on the boundary of the 3D bulk, should be enough to train the neural network. It would also be interesting to extend the current framework for time-dependent physical systems, for example, viscoelastic materials or elasto-plastic materials. In that case, the constitutive relations are path-dependent, and the deformation in time of the material is used as training data. In such cases, it may be possible to treat the forward propagation as a standalone operator and use the adjoint state method to derive the gradient of the operator. In this way, the neural network approximation to unknown functions is decoupled from the sophisticated numerical simulations. Besides constitutive relations, other coarse-grained models in damage mechanics and structural internal damping will be explored in the future.


\section*{Acknowledgements}
This work is supported by the Applied Mathematics Program within the Department of Energy (DOE) Office of Advanced Scientific Computing Research (ASCR), through the Collaboratory on Mathematics and Physics-Informed Learning Machines for Multiscale and Multiphysics Problems Research Center (DE-SC0019453). Daniel Z.\ Huang thanks Prof.\ Peter Pinsky for providing the model problems in our numerical experiments. Kailai Xu also thanks the Stanford Graduate Fellowship in Science \& Engineering and the 2018 Schlumberger Innovation Fellowship for partial financial support.

\section*{Appendix}

\added{The library \texttt{ADCME} is built on \texttt{TensorFlow}. Through operator overloading and multiple dispatch features of \texttt{Julia}, users can almost write \texttt{Julia} style syntax for creating computational graph in \texttt{TensorFlow}, which is used for automatic differentiation. Additional, missing features in \texttt{TensorFlow} that are essential for scientific computing such as sparse solvers are implemented though \texttt{C++} kernels and ``plugged'' in as custom operators. Numerical operations that do not require gradient back-propagation are executed in \texttt{Julia}, which are more efficient as opposed to that implemented in pure \texttt{Python}.}

In this appendix, we explain how to solve a simple inverse modeling problem in \texttt{ADCME} and how automatic differentiation works in the \texttt{TensorFlow} backend. We consider the following Poisson equation where we want to estimate the unknown scalar $b$ in 
\begin{equation}\label{equ:model problem 2}
    -b u''(x) + u(x) = f(x), \quad x\in [0,1], \quad u(0)=u(1)=0
\end{equation}
where
\begin{equation}
    f(x) = 8+4x-4x^2
\end{equation}
Assume that the true parameter $b^*=1$ and we have observed the corresponding solution $u(x)$ at $x=0.5$, i.e., $u(0.5)=1$.

\begin{figure}[hbt]
  \includegraphics[width=1.0\textwidth]{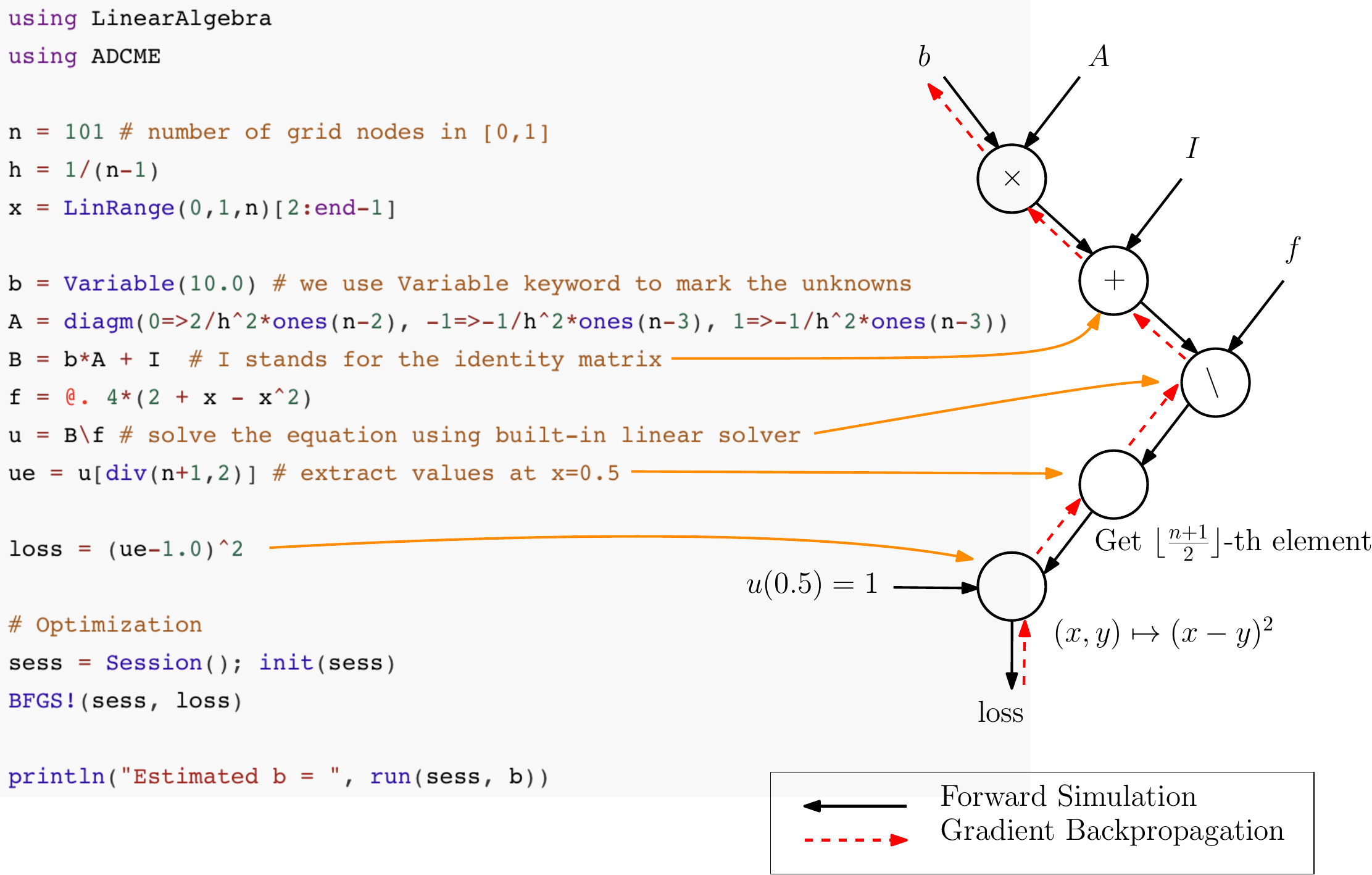}
  \caption{Illustration of solving a simple inverse modeling problem in \texttt{ADCME}.}
  \label{fig:code2}
\end{figure}

\texttt{TensorFlow} represents data by Tensors, which are multidimensional arrays with extra traits such as data types, shapes, and whether they are trainable or not. In principle, all the data to be processed by the \texttt{TensorFlow} backend must be converted to Tensors; however, we have overloaded \texttt{Julia} operators for convenience so that \texttt{Julia} arrays are compatible with Tensors, e.g., we can add a Tensor and a Julia array without first converting the latter to a Tensor. A Variable is a special Tensor marked as ``trainable''. During optimization, \texttt{TensorFlow} looks for those trainable Variables and computes the gradients of the objective function with respect to them. The values of Variables are updated during the optimization process. 

For example, in the code snippet in \cref{fig:code2}, since $b$ is unknown, we create a Variable with the initial guess $b=10$ 
\begin{center}
    \texttt{b = Variable(10.0)}
\end{center}
other data such as the coefficient matrix $A$, the identity matrix $I$, the source term $f$, and the data $u(0.5)=1$ are all Tensors~(not trainable), but programmatically we can represent them with \texttt{Julia} arrays thanks to operator overloading.

The key to automatic differentiation is the construction of a computational graph. A computational graph represents each \texttt{TensorFlow} operation (such as $\times$, $+$, matrix solve, indexing, etc.) as a node and each Tensor (either input or intermediate result) as an edge. Each node~(operation) takes zero or several Tensors as inputs and outputs a Tensor. This convention (nodes represent operations while edges represent Tensors) ensures that no Tensor is the output of multiple operations, and operations can take multiple Tensors as inputs. For example, in \Cref{fig:code2} the multiplication operator $\times$ takes two inputs, $b$ and $A$, and outputs another Tensor $bA$. 

The computational graph keeps track of the computation dependencies that are used for ``back-propagating'' gradients. To understand the back-propagation process, we look at a single operator $y = F(x)$ where $x$ is the input Tensor while $y$ is the output Tensor. The scalar loss $l$ depends on $y$ and thus on $x$, i.e. $l = L(y) = L(F(x))$. To compute the gradients $\frac{\partial L(F(x))}{\partial x}$, we have
\begin{equation}\label{equ:prop}
    \frac{\partial L(F(x))}{\partial x} = \frac{\partial L(y)}{\partial y} \frac{\partial F\added{(x)}}{\partial x}
\end{equation}
the quantity $\frac{\partial L(y)}{\partial y}$ is ``back-propagated'' from the loss function, and the operator $F$ computes gradients $\frac{\partial L(F(x))}{\partial x} = \frac{\partial L(y)}{\partial y} \frac{\partial F(x)}{\partial x}$ and passes it to the next operator. By next operator, we mean that $x$ may be the output of another operator $G$; this operator $G$ takes $\frac{\partial L(F(x))}{\partial x}$ as ``back-propagated'' gradients from the loss function and passes the resulting gradients likewise. Therefore, by specifying the propagating rule \cref{equ:prop} for each operator, we can automate gradient computations. 

Finally, we explain how to solve the model problem \cref{equ:model problem 2} in \texttt{ADCME}. We discretize the system by the finite difference method on a uniform grid $0 = x_1 < x_2 < \cdots < x_{n+1}  = 1$~($n=100$) and denote the approximated nodal values at $x_i$ by $u_i$. The discretized system corresponding to \cref{equ:model problem 2} is (after taking into account the boundary condition $u_1=u_{n+1}=0$)
\begin{equation}\label{equ:discretized}
    (bA+I)\bu = \mathbf{f}
\end{equation}
where
\begin{equation}
    A = \begin{bmatrix}
        \frac{2}{h^2} & -\frac{1}{h^2} & \dots & 0\\
         -\frac{1}{h^2} & \frac{2}{h^2} & \dots & 0\\
         \dots \\
         0 & 0 & \dots & \frac{2}{h^2}
    \end{bmatrix}, \quad \mathbf{u} = \begin{bmatrix}
        u_2\\
        u_3\\
        \vdots\\
        u_{n}
    \end{bmatrix}, \quad \mathbf{f} = \begin{bmatrix}
        f(x_2)\\
        f(x_3)\\
        \vdots\\
        f(x_{n})
    \end{bmatrix}
\end{equation}
where $I$ is the identity matrix. The idea is to solve \Cref{equ:discretized} given the \texttt{Variable} $b$ to obtain $\bu$, and minimize the discrepancy between the observation $u(0.5)=1$ and $u_{\lfloor \frac{n+1}{2} \rfloor}$. The following code shows the process of solving the inverse problem,

\begin{enumerate}
    \item The first two lines in \cref{fig:code2} load necessary packages
\begin{lstlisting}
using LinearAlgebra
using ADCME    
\end{lstlisting}
\item We split the interval $[0,1]$ into $100$ equal length subintervals
\begin{lstlisting}
n = 101
h = 1/(n-1)
x = LinRange(0,1,n)[2:end-1]
\end{lstlisting}
\item Since $b$ is unknown and needs to be updated during optimization, we mark it as trainable using the \texttt{Variable} keyword
\begin{lstlisting}
b = Variable(10.0)
\end{lstlisting}
\item Solve the linear system \Cref{equ:discretized} and extract the value at $x=0.5$
\begin{lstlisting}
A = diagm(0=>2/h^2*ones(n-2), -1=>-1/h^2*ones(n-3), 1=>-1/h^2*ones(n-3)) 
B = b*A + I  
f = @. 4*(2 + x - x^2) 
u = B\f 
ue = u[div(n+1,2)] 
\end{lstlisting}
here \texttt{I} stands for the identity matrix and \texttt{@.} denotes element-wise operation. They are \texttt{Julia}-style syntax but are also compatible with tensors by overloading. 
\item Formulate the loss function
\begin{lstlisting}
loss = (ue-1.0)^2 
\end{lstlisting}
\item Create and initialize a \texttt{TensorFlow} session, which analyzes the computational graph and initializes the tensor values. 
\begin{lstlisting}
sess = Session(); init(sess) 
\end{lstlisting}
\item Finally, we trigger the optimization by invoking \texttt{BFGS!}, which wraps the \texttt{L-BFGS-B} algorithm
\begin{lstlisting}
BFGS!(sess, loss)
\end{lstlisting}
since the computational graph contains all the information, this process, including gradient computation, variable update, etc., has been fairly automated. 
\end{enumerate}

\bibliographystyle{unsrt}
\bibliography{PMUQ}

\end{document}